\documentclass{amsart}
\title{Extender Based Radin Forcing}
\author{Carmi Merimovich}
\address{
School of Mathematical Sciences
\\
Tel-Aviv University
\\
Tel-Aviv 69978
\\
ISRAEL
}
\email{carmi\_m@cet.ac.il}
\date{October 10, 1998}

\thanks{
This work is a part of research which, hopefully, will become the author's
Ph.D. thesis. It was done at Tel-Aviv university under the supervision
of M. Gitik. The author thanks M. Gitik for his help with this work,
in other works and just in general.
}


\def\resurrect{0}

\usepackage{amssymb}

\usepackage{pb-diagram}

\theoremstyle{plain}
\newtheorem{theorem}{Theorem}[section]
\newtheorem{lemma}[theorem]{Lemma}
\newtheorem{proposition}[theorem]{Proposition}
\newtheorem{claim}[theorem]{Claim}
\newtheorem{corollary}[theorem]{Corollary}

\theoremstyle{definition}
\newtheorem{definition}[theorem]{Definition}

\theoremstyle{remark}

\newtheorem{note}[theorem]{Note}

\numberwithin{equation}{theorem}

\newcommand{\cA}{{\mathcal{A}}}

\newcommand{\ga}{\alpha}
\newcommand{\gb}{\beta}
\newcommand{\gc}{\chi}

\newcommand{\gee}{\epsilon}

\newcommand{\gga}{\gamma}

\newcommand{\gk}{\kappa}
\newcommand{\gl}{\lambda}
\newcommand{\gm}{\mu}
\newcommand{\gn}{\nu}

\newcommand{\gp}{\pi}

\newcommand{\gr}{\rho}
\newcommand{\gs}{\sigma}
\newcommand{\gt}{\tau}

\newcommand{\gw}{\omega}
\newcommand{\gx}{\xi}

\newcommand{\gz}{\zeta}

\newcommand{\func}{\mathord{:}}
\newcommand{\VN}[1]{\widehat{#1}}
\newcommand{\GN}[1]{\widetilde{#1}}
\newcommand{\satisfies}{\vDash}
\newcommand{\union}{\cup}
\newcommand{\bigunion}{\bigcup}
\newcommand{\intersect}{\cap}
\newcommand{\bigintersect}{\bigcap}
\newcommand{\forces}{\mathrel\Vdash}
\newcommand{\incompatible}{\perp}
\newcommand{\compatible}{\parallel}
\newcommand{\decides}{\mathrel\Vert}

\newcommand{\subelem}{\prec}

\newcommand{\append}{\mathop{{}^\frown}}
\newcommand{\restricted}{\mathord{\restriction}}
\newcommand{\upto}{\mathord{<}}
\newcommand{\downto}{\mathord{>}}
\newcommand{\power}[1]{\lvert#1\rvert}
\newcommand{\Pset}{{\mathcal{P}}}

\newcommand{\ordered}[1]{\ensuremath{\langle #1 \rangle}}
\newcommand{\bordered}[1]{\ensuremath{\big\langle #1 \big\rangle}}

\newcommand{\set}[1]{\ensuremath{\left\{ #1 \right\}}}
\newcommand{\setof}[2]{\ensuremath{\left\{ #1 \mid #2 \right\}}}
\newcommand{\ordof}[2]{\ensuremath{\ordered{ #1 \mid #2 }}}
\newcommand{\formula}[1]{{}^{\ulcorner} #1 {}^{\urcorner}}

\DeclareMathOperator{\len}{l}
\DeclareMathOperator{\otp}{otp}
\DeclareMathOperator{\crit}{crit}
\DeclareMathOperator{\dom}{dom}
\DeclareMathOperator{\cf}{cf}
\DeclareMathOperator{\supp}{supp}
\DeclareMathOperator{\Ult}{Ult}
\DeclareMathOperator{\Suc}{Suc}
\DeclareMathOperator{\Lev}{Lev}
\DeclareMathOperator{\ran}{ran}
\DeclareMathOperator*{\dintersect}{\triangle}
\newcommand{\dsintersect}{\sideset{}{^0}\dintersect}

\DeclareMathOperator{\mc}{mc}

\newcommand{\gas}{{\ensuremath{\Bar{\ga}}\/}}
\newcommand{\gbs}{{\ensuremath{\Bar{\gb}}\/}}
\newcommand{\ggs}{{\ensuremath{\Bar{\gga}}\/}}
\newcommand{\ges}{{\ensuremath{\Bar{\gee}}\/}}
\newcommand{\gnv}{{\ensuremath{\Vec{\gn}}\/}}
\newcommand{\gns}{{\ensuremath{\Bar{\gn}}\/}}
\newcommand{\gms}{{\ensuremath{\Bar{\gm}}\/}}

\newcommand{\Es}{{\ensuremath{\bar{E}}\/}}
\newcommand{\Pe}{{\ensuremath{P_{\ges}\/}}}
\newcommand{\PE}{{\ensuremath{P_{\Es}\/}}}

\newcommand{\Cond}[2]{\ordered{\Vec{#1},\Vec{#2}}}

\begin{document}
\begin{abstract}
We define extender sequences,  generalizing
measure sequences of Radin forcing.

Using the extender sequences, we show how to combine the 
Gitik-Magidor forcing
for adding many Prikry sequenes with Radin forcing.

We show that this forcing satisfies a Prikry like condition
destroys no cardinals, and has a kind of properness.

Depending on the large cardinals we start with
this forcing can blow the power of a cardinal together
with changing its' cofinality to a prescirbe value.
It can even blow the power of a cardinal while keeping it regular
or measurable.
\end{abstract}
\maketitle
\section{Introduction}
We give some background on previous work
relating directly to the present work.
The first forcing which changed
the cofinality of a cardinal without changing the cardinal structure
was Prikry forcing \cite{Prikry}. In this forcing a measurable cardinal, $\gk$,
was `invested' in order to get $\cf(\gk)=\gw$ without collapsing
any cardinal.
Developing that idea, Magidor \cite{Menachem:cf}
used a coherent sequence
of measures of length $\gl < \gk$ in 
order to get $\cf(\gk)=\gl$ without collapsing any cardinals.
In \cite{Radin} Radin, introducing the notion of measure sequence, showed
that it is useful to continue the coherent sequence to $\gl > \gk$.
For example, $\gk$ remains regular when $\gl = \gk^+$.
In general
the longer the measure sequence the more resemblance there is
between $\gk$ in the generic extension and the ground model.

As is well known, and unlike regular cardinals, 
blowing the power of a singular cardinal is not
an easy task. A natural approach to try was to blow the power
of a cardinal while it was regular and after that  make it singular
by one of the above methods.
A crucial idea of Gitik and Magidor \cite{Moti} was to combine
the power set blowing and the cofinality change in one forcing.
They introduced a forcing notion which added many Prikry
sequences at once and still collapsed no cardinals. 
The `investment' they needed for this was an extender
of length which is the size of power they wanted.
Building on the idea of Gitik and Magidor, Segal \cite{Miri} 
implemented the idea of adding many sequences
to Magidor forcing. So by investing a coherent sequence of extenders
of length $\gl < \gk$ she was able to get a singular cardinal of
cofinality $\gl$ together with power as large as the length of the
extenders in question.
Our work also builds on the idea of Gitik and Magidor.
However, we implement the idea of adding many sequences to Radin forcing.
So we introduce the notion of extender sequence and show that it makes
sense to deal with  quite long extender sequences. As in Radin forcing,
for long enough sequences we are left with $\gk$ which is regular and
even measurable. The power size will be the length 
of the extenders we  start with.

The structure of this work is as follows. In section \ref{ExtenderSequences}
we define extender sequences.
In section \ref{RadinForcing} we define Radin forcing. The definition
is not the usual one and is used to introduce the idea
used in section \ref{PEForcing}.
In section \ref{PEForcing} we define $\PE$, the forcing notion which is the
purpose of this work.
In section \ref{BasicProperties} we show the chain-conditions satisfied
by $\PE$ and how `locally' it resembles Radin forcing. We also show here
that there are many new subsets in the generic extension.
In section \ref{HomogenDense} we investigate the structure of dense open
subsets of $\PE$. We show that they satisfy
a strong homogeneity property. 
In section \ref{Prikry'sCondition} we prove Prikry's condition for $\PE$.
The proof is a simple corollary of the
strong homogeneity of dense open subsets.
In section \ref{Properness} 
we show that $\PE$ satisfies
a kind of properness. 
In section \ref{Cardinals} we combine the machinery developed so far
in order to show that no cardinals are collapsed.
In section \ref{PropertiesK} we show how the length of the extender
sequences affect the properties of $\gk$.
Section \ref{TheTheorem} summarizes
what the forcing
$\PE$ does.
In section \ref{ByIteration} we have a result concerning $\PE$ when
$\len(\Es)=1$. We show that 
there is, in $V$, a generic filter
over an elementary submodel in an $\gw$-iterate of $V$. We were not able
to prove something equivalent (or weaker) for the general case.
Section \ref{ConcludingRemarks} contains a list of missing or
unknown points to check. The last point in this list is in preparation.

Our notation is standard. We assume fluency with forcing and
extenders. Some basic properties of Radin forcing are taken for
granted.
\section{Extender sequences} \label{ExtenderSequences}
\subsection{Constructing from elementary embedding}
Suppose we have an elementary embedding $j\func V \to M \supset V_\gl$, $\crit(j)=\gk$.
The value of $\gl$ is determined later, according to the different
applications we have.

Construct from $j$ a nice extender like in \cite{Moti}:
\begin{align*}
E(0) = \ordered{\ordof{E_\ga(0)}{\ga \in \cA},
                \ordof{\gp_{\gb,\ga}}{\gb \geq \ga \  \ga,\gb \in \cA}}.
\end{align*}
We recall the properties of this extender:
\begin{enumerate}
\item $\cA \subseteq \power{V_\gl}\setminus\gk$,
\item $\power{\cA} = \power{V_\gl}$,
\item $\cA$ is $\gk^+$-directed,
\item $\gk$ is minimal in $\cA$ and we write $\gp_{\ga,0}$ instead
        of $\gp_{\ga,\gk}$,
\item $\forall \ga,\gb\in \cA$ 
        $\gn^0 = \gp_{\ga,0}(\gn) = \gp_{\gb,0}(\gn)$,
\item $\forall \ga,\gb \in \cA$ 
        $\gp_{\gb,0}(\gn) = \gp_{\ga,0}(
                                \gp_{\gb,\ga}(\gn))$,
\item $\forall \ga,\gb,\gga \in \cA$ $\exists A \in E_\gga(0)$
        $\forall \gn \in A$
        $\gp_{\gga,\ga}(\gn) = \gp_{\gb,\ga}(
                                \gp_{\gga,\gb}(\gn))$.
\end{enumerate}
If, for example, we need $\power{E(0)} = \gk^{+3}$ then, under $\text{GCH}$,
we require $\gl = \gk+3$.
A typical large set in this extender concentrates on singletons.

If $j$ is not sufficiently closed , then $E(0) \notin M$ and the construction
stops. We set
\begin{align*}
\forall \ga \in \cA \ \Es_\ga = \ordered{\ga,E(0)}.
\end{align*}
We say that $\Es_\ga$ is an extender sequence of length $1$.
($\len(\Es_\ga)=1$)

If, on the other hand, $E(0) \in M$ we can construct for each $\ga \in \dom E(0)$
the following ultrafilter
\begin{align*}
A \in E_{\ordered{\ga,E(0)}}(1) \iff \ordered{\ga, E(0)} \in j(A).
\end{align*}
Such an $A$ concentrates on elements of the form $\ordered{\gx, e(0)}$
where $e(0)$ is an extender on $\gx^0$
and $\gx \in \dom e(0)$. Note that $e(0)$  concentrates on singletons
below $\gx^0$. If, for example, $\power{E(0)} = \gk^{+3}$ then on a large
set we have $\power{e(0)} = (\gx^{0})^{+3}$.

We define $\gp_{\ordered{\gb,E(0)},\ordered{\ga,E(0)}}$ as
\begin{align*}
\gp_{\ordered{\gb,E(0)},\ordered{\ga,E(0)}}(\ordered{\gx, e(0)}) = 
        \ordered{\gp_{\gb,\ga}(\gx), e(0)}.
\end{align*}
From this definition we get
\begin{align*}
j(\gp_{\ordered{\gb,E(0)},\ordered{\ga,E(0)}})(\ordered{\gb, E(0)}) = 
        \ordered{\ga, E(0)}.
\end{align*}
Hence we have here an extender
\begin{align*}
E(1) = \ordered{
                \ordof{E_{\ordered{\ga,E(0)}}(1)}{\ga \in {\cA}},
                \ordof{\gp_{\ordered{\gb,E(0)},\ordered{\ga,E(0)}}}
                        {\gb \geq \ga \  \ga,\gb \in \cA}
        }.
\end{align*}
Note that the difference between $\gp_{\gb,\ga}$ and
$\gp_{\ordered{\gb,E(0)},\ordered{\ga,E(0)}}$ is quite superficial.
We can define $\gp_{\ordered{\gb,E(0)},\ordered{\ga,E(0)}}$
in a uniform way for both extenders. Just project the first element of
the argument using $\gp_{\gb,\ga}$.

If $\ordered{E(0),E(1)} \notin M$ then the construction stops.
In this case we set
\begin{align*}
\forall \ga \in \cA \ \Es_\ga = \ordered{\ga, E(0), E(1)}.
\end{align*}
We say that $\Es_\ga$ is an extender sequence of length $2$.
($\len(\Es_\ga)=2$)

If $\ordered{E(0),E(1)} \in M$ then we construct the extender
$E(2)$ in the same way as we constructed $E(1)$ from $E(0)$.

The above private case being worked out we continue with the general case.
Assume we have constructed
\begin{align*}
\ordof  {E(\gt')}
        {\gt' < \gt}.
\end{align*}
If $\ordof {E(\gt')} {\gt' < \gt} \notin M$ then the construction stops
here. We set
\begin{align*}
\forall \ga \in \cA \ \Es_\ga = \ordof  {\ga, E(\gt')}
                        {\gt' < \gt},
\end{align*}
and we say that $\Es_\ga$ is
an extender sequence of length $\gt$. ($\len(\Es_\ga)= \gt$)

If, on the other hand, $\ordof {E(\gt')} {\gt' < \gt} \in M$ then 
we construct
\begin{multline*}
A \in E_{\ordof{\ga, E(0),\dotsc, E(\gt'),\dotsc} {\gt' < \gt}}(\gt) \iff
\\
        \ordof{\ga, E(0), \dotsc, E(\gt'), \dotsc} {\gt' < \gt}
                                         \in j(A).
\end{multline*}
Defining $\gp_{\ordof{\gb, E(0),\dotsc, E(\gt'),\dotsc}{\gt' < \gt},
\ordof{\ga, E(0),\dotsc, E(\gt'),\dotsc,}{\gt' < \gt}}$ using the first
coordinate as before gives the needed projection.

We are quite casual in writing the indices of the projections and
ultrafilters.
By this we mean that we  sometimes write $\pi_{\gb,\ga}$ when we
should have written
$\gp_{\ordof{\gb, E(0),\dotsc, E(\gt'),\dotsc}{\gt' < \gt},
\ordof{\ga, E(0),\dotsc, E(\gt'),\dotsc,}{\gt' < \gt}}$
and $E_\ga(\gt)$ when we should have written
$E_{\ordof{\ga, E(0),\dotsc, E(\gt'),\dotsc,}{\gt' < \gt}}(\gt)$.

With this abuse of notation the projection we just defined satisfies
\begin{multline*}
j(\gp_{\gb,\ga})
        (\ordof{\gb, E(0), \dotsc, E(\gt'), \dotsc} {\gt' < \gt})=
\\
                \ordof{\ga, E(0), \dotsc, E(\gt'), \dotsc} {\gt' < \gt},
\end{multline*}
and we have the extender
\begin{align*}
E(\gt) = \ordered{\ordof{E_\ga(\gt)}{\ga \in {\cA}},
                \ordof{\gp_{\gb,\ga}}{\gb \geq \ga \  \ga,\gb \in \cA}}.
\end{align*}

We let the construction run until it stops due to the extender
sequence not being in $M$.

\begin{definition}
We call $\gms$ an extender sequence if there is an elementary embedding
$j\func V \to M$ such that $\gns$ is an extender sequence generated 
as above and $\gms = \gns \restricted \gt$ for $\gt \leq \len(\gns)$.
$\gk(\gms)$ is the ordinal at the beginning of
the sequence. (i.e. $\gk(\bar{E}_\ga)=\ga$). $\gk^0(\gms)$ is
$(\gk(\gms))^0$. (i.e. $\gk^0(\bar{E}_\ga)=\gk$).
\end{definition}
That is, we do not have to construct the extender sequence until it is
not in $M$. We can stop anywhere on the way.
\begin{definition}
A sequence of extender sequences $\ordered{\gms_1, \dotsc, \gms_n}$ is
called $^0$-increasing if $\gk^0(\gms_1) < \dotsb <\gk^0(\gms_n)$.
\end{definition}
\begin{definition}
Let $\ordered{\gms_1, \dotsc, \gms_n}$ be $^0$-increasing.
An extender sequences $\gms$ is called permitted to 
$\ordered{\gms_1, \dotsc, \gms_n}$ if
         $\gk^0(\gms_n) < \gk^0(\gms)$.
\end{definition}
%
% When a Set is in an extender sequence (definition)
%
\begin{definition}
We say $A \in \bar{E}_\ga$ if $\forall \gx < \len(\bar{E}_\ga) \ A\in E_\ga(\gx)$.
\end{definition}
\begin{definition}
$\Es =\ordof{\Es_\ga}{\ga \in \cA}$ is an extender sequence system
if there is an elementary embedding $j\func V \to M$ such that all
$\Es_\ga$ are extender sequences generated from $j$ as prescribed above
and $\forall \ga,\gb\in{\cA} \ \len(\Es_\ga)=\len(\Es_\gb)$. This
common length is called the length of the system, $\len(\Es)$. We
write $\Es(\gms)$ for the 
extender sequence system to which $\gms$ belongs (i.e. $\Es(\Es_\ga)=\Es$).
\end{definition}
The generalization of the measure on the $\ga$ coordinate in
Gitik-Magidor forcing \cite{Moti} is $\Es_\ga$.
%
% E-tree subsection
%
\subsection{$\bar{E}_\ga$-tree}
%
% E-tree defintion
%
\begin{definition}
A tree $T$ is an $\bar{E}_\ga$-tree if its' elements are
of the form 
\begin{equation*}
\ordered{\ordered{\gms_1,\linebreak[0] \dotsc, \gms_n},S}
\end{equation*}
where
\begin{enumerate}
\item Set $\dom T = \setof{\ordered{\gms_1,\dotsc,\gms_n}}
        {\ordered{\ordered{\gms_1,\dotsc,\gms_n},S} \in T}$. Then the function
        \begin{align*}
        \ordered{\gms_1,\dotsc,\gms_n} \mapsto \ordered{\ordered{\gms_1,\dotsc,\gms_n},S}
        \end{align*}
        from $\dom T$ to $T$ is $1-1$ and onto,

\item $t \in \Lev_n(\dom T) \implies \power{t}=n+1$,

\item $\ordered{\gms_1,\dotsc,\gms_n}$ are $^0$-increasing extender sequences,

\item $\Lev_0(\dom T)\in \bar{E}_\ga$ and for each $t\in \dom T$
      $\Suc_{\dom T}(t)\in \bar{E}_\ga$,

\item $S$ is a $\gms_n$-tree. When $\len(\gms_n)=0$ we 
      set $S=\emptyset$.
      \par\noindent Note that this clause is recursive.
\end{enumerate}
\end{definition}
\begin{note}
Later on, we abuse notation and use $T$ instead of $\dom T$.
i.e. $\Suc_T(t)$ instead of $\Suc_{\dom T}(t)$.
\end{note}
%
% Define sub-trees
%
\begin{definition}
Assume $T$ is a $\bar{E}_\ga$-tree and $t\in T$, then:
\begin{enumerate}
\item $T_t=\setof{\ordered{s, S}} {\ordered{t\append s, S} \in T}$.
\item $T_t(\gms)$ is the tree, $S$, satisfying
      $\ordered{\gms, S}\in \Suc_T(t)$.
\end{enumerate}
\end{definition}
%
% Define order on Prikry trees
%
\begin{definition}
Let $T$, $S$ be $\bar{E}_\ga$-trees, where $\len(\bar{E})=1$. We say that
$T \leq S$ if
\begin{enumerate}
\item $\Lev_0(T) \subseteq \Lev_0(S)$,
\item $\forall t \in T \Suc_T(t)\subseteq \Suc_S(t)$.
\end{enumerate}
\end{definition}
%
% Define order on Magidor-Radin trees
%
\begin{definition}
Let $T$, $S$ be $\bar{E}_\ga$-trees. We say that $T \leq S$ if
\begin{enumerate}
\item $\Lev_0(T) \subseteq \Lev_0(S)$,
\item $\forall t \in T \Suc_T(t)\subseteq \Suc_S(t)$,
\item $\forall \ordered{\gms}\in\Lev_0(T) \ T(\gms)\leq S(\gms)$,
\item $\forall t \in T \ \forall \ordered{\gms} \in T_t \ 
                        T_t(\gms)\leq S_t(\gms)$.
\end{enumerate}
Note that the last 2 conditions are recursive.
\end{definition}
\begin{definition}
Let $S$ be $\bar{E}_\ga$-tree and $\gb>\ga$. Define
$T = \pi^{-1}_{\gb,\ga}(S)$ by
\begin{enumerate}
\item $\dom T = \gp^{-1}_{\gb,\ga} (\dom S)$,
\item $T_{\ordered{\gms_1,\dotsc,\gms_{n-1}}}(\gms_n) = 
        \gp^{-1}_{\gms_n, \gp_{\gb,\ga}(\gms_n)}
                S_{\ordered{\gp_{\gb,\ga}(\gms_1),\dotsc,
                          \gp_{\gb,\ga}(\gms_{n-1})}}(\gp_{\gb,\ga}(\gms_{n}))$.
\end{enumerate}
\end{definition}
\begin{definition}
Let $T$, $S$ be $\bar{E}_\gb$, $\bar{E}_\ga$-trees respectively, where
$\gb \geq \ga$. We say that $T \leq S$ if
\begin{enumerate}
\item $T \leq \pi^{-1}_{\gb,\ga} (S)$.
\end{enumerate}
\end{definition}
%
% Define Diagonal Intersection
%
\begin{definition}
Assume we have $A_{\ordered{\gns, R}}$ where $\gns$ is an extender sequence
such that each element in $A_{\ordered{\gns,R}}$
is of the form $\ordered{\gms, S}$ where $\gms$ is an extender
sequence and $S$ is a tree.
(in this work we always have $S$ is $\gms$-tree).
We define $\dintersect^0_{\ordered{\gns,R}} A_{\ordered{\gns,R}}$ as
\begin{align*}
\ordered{\gms,S} \in \sideset{}{^0}\dintersect_{\ordered{\gns,R}} A_{\ordered{\gns,R}}
        \iff
        \forall \ordered{\gns,R} \, \gk(\gns) < \gk^0(\gms) 
        \rightarrow
         \ordered{\gms,S_{\ordered{\gns, R}}} \in A_{\ordered{\gns,R}}.
\end{align*}
\end{definition}
%
% R A D I N    F o r c i n g
%
\section{Radin Forcing} \label{RadinForcing}
The main aims of this section are \ref{fill-missing}, \ref{skeleton}.
As the simplest way we found to formulate them was with Radin forcing
\cite{Radin,Mitchell,Woodin} we took the opportunity
to depart from  usual formulation in order to introduce 
ideas that we use in the extender based forcing later.

The main point is that possible extensions of a condition are stored in 
$\Es_\ga$-tree 
and not in a set. The $\ga$ is be fixed, so practically we deal here
with a measure sequence and not an extender sequence.

\begin{definition}
A condition in $R_\ga$ is of the form
\begin{align*}
\ordered{\ordered{\gms_n,s^n},S^n, \dotsc, \ordered{\gms_{1},s^1},S^{1}, 
        \ordered{\gms_0,s^0}, S^0}
\end{align*}
where
\begin{enumerate}
\item $\gms_0 = \Es_\ga$,
\item $\forall i\leq n$ $S^i$ is $\gms_i$-tree,
\item $\forall i\leq n$ 
        $s^i \in V_{\gk^0(\gms_i)}$ is an extender sequence.
\end{enumerate}
\end{definition}
\begin{definition}
Let $p,q \in R_\ga$. We say that $p$ is Prikry extension of $q$
($p \leq^* q$ or $p \leq^0 q$) if $p,q$ are of the form
\begin{align*}
& p = \ordered{\ordered{\gms_n,s^n},S^n, \dotsc, \ordered{\gms_{1},s^1},S^{1},
                \ordered{\Es_\ga,s^0}, S^0},
\\
& q = \ordered{\ordered{\gms_n,t^n},T^n, \dotsc, \ordered{\gms_{1},t^1},T^{1},
                \ordered{\Es_\ga,t^0},T^0},
\end{align*}
and
\begin{itemize}
\item $\forall i \leq n$ $S^i \leq T^i$,
\item $\forall i \leq n$ $s^i = t^i$.
\end{itemize}
\end{definition}
\begin{definition}
Let
$p = \ordered{\ordered{\gms_n,s^n},S^n, \dotsc, \ordered{\gms_{1},s^1},S^{1},
                \ordered{\gms_0,s^0}, S^0}$
where $\gms_0=\Es_\ga$. Let $\ordered{\gns}\in S^i$. We define 
$(p)_{\ordered{\gns}}$ to be
\begin{align*}
%%\begin{split}
(p)_{\ordered{\gns}} = 
        \ordered{
        \ordered{\gms_n,s^n},S^n, \dotsc, 
        \ordered{\gms_{i+1},s^{i+1}},S^{i+1},
\\
        \ordered{\gns, s^i}, S^i(\gns),
        \ordered{\gms_i, \gns}, S^i_{\ordered{\gns}},
\\
        \ordered{\gms_{i-1},s^{i-1}},S^{i-1},
        \dotsc,
        \ordered{\gms_0,s^0}, S^0
        }.
%%\end{split}
\end{align*}
\end{definition}
Note the degenerate case in this definition when $\len(\gns)=0$. In this
case $S^i(\gns)=\emptyset$.
\begin{definition}
Let $p,q \in R_\ga$ where 
\begin{align*}
q = \ordered{\ordered{\gms_n,s^n},S^n, \dotsc, \ordered{\gms_{1},s^1},S^{1},
                \ordered{\gms_0,s^0}, S^0}.
\end{align*}
We say that $p$ is 1-point extension of $q$
($p \leq^1 q$) if there is $\ordered{\gns} \in S^i$ such that
$p \leq^* (q)_{\ordered{\gns}}$.
\end{definition}
\begin{definition}
Let $p,q \in R_\ga$. We say that $p$ is $n$-point extension of $q$
($p \leq^n q$) if there are $p^n, \dotsc, p^0$ such that
\begin{align*}
p = p^n \leq^1 \dotsc \leq^1 p^0 = q.
\end{align*}
\end{definition}
\begin{definition}
Let $p,q \in R_\ga$. We say that $p$ is an extension of $q$
($p \leq q$) if there is $n$ such that $p \leq^n q$.
\end{definition}
%
% Tree completion lemma
%
\ref{fill-missing} is needed in the proof of \ref{PrikryCondition}.
Very loosly speaking \ref{fill-missing} means that if ``something''
happens on a measure one set for one of the measures, that
``something'' is happening on a measure one set for all the measures.

\ref{fill-missing} is proved by induction and \ref{fill-missing2} is
the first case of the induction.
\begin{lemma} \label{fill-missing2}
Suppose $\len(\Es_\ga)=2$ and $T$ is a tree such that
$\Lev_0(T)\in E_\ga(i)$
for $i<2$, and $\forall \ordered{\gns}\in T$
$T_{\ordered{\gns}}$ is an $\Es_{\ga}$-tree.
Then there is an $\Es_{\ga}$-tree, $T^*$, satisfying
\begin{enumerate}
\item $\forall \ordered{\gns}\in T \intersect T^*$  
      $T^*_{\ordered{\gns}} \leq T_{\ordered{\gns}}$,
      $T^*(\gns) \leq T(\gns)$,
\item If 
        $p
                \leq
              \bordered{\ordered{\Es_\ga,\ordered{}}, T^*}$ 
        then there is 
        $\ordered{\gms}\in T^* \intersect T$
       such that
       \begin{align*}
      p
        \compatible
      \big(
        \bordered{\ordered{\Es_\ga,\ordered{}}, T^*}
                \big)_{\ordered{\gms}}.
       \end{align*}
\end{enumerate}
\end{lemma}
\begin{proof}
There are 2 cases which to deal with:
\begin{itemize}
\item $A_0 = \Lev_0(T) \in E_\ga(0)$:
If $A_0 \in E_\ga( 1)$ we set $T^*=T$ and the 
proof is finished. So suppose
$A_0 \notin E_\ga( 1)$. We would like to
build $A_1 \in E_\ga( 1)$. Set
\begin{align*}
& \forall \ordered{\gms_0} \in T \ A_{\ordered{\gms_0}, 1}=
        \setof{\ordered{\gms_1,S} \in T_{\ordered{\gms_0}}}
              {A_0 \intersect \gk^0(\gms_1) \in
                                \gms_1(0)_{\gk(\gms_1)}}.
\end{align*}
As $A_0 \in E_\ga(0)$ and $\Suc_T(\ordered{\gms_0}) \in E_\ga(1)$
we get that $A_{\ordered{\gms_0},1} \in E_\ga(1)$.
Let
\begin{align*}
& A_1 = \dsintersect_{\ordered{\gms_0}\in T} A_{\ordered{\gms_0},1}.
\end{align*}
We can construct now $T^*$:
\begin{align*}
&\Lev_0(T^*)= A_0 \union A_1,
\\
&\forall\gms_0 \in A_0 \ T^*_{\ordered{\gms_0}}= T_{\ordered{\gms_0}},
\\
&\forall\gms_1 \in A_1\ T^*_{\ordered{\gms_1}}=
\bigintersect_{\ordered{\gms_0,\gms_1}\in T}T_{\ordered{\gms_0,\gms_1}}.
\end{align*}
\item $A_1 = \Lev_0(T)\in E_\ga( 1)$:
If $A_1 \in
E_\ga( 0)$ we set $T^*=T$ and finish the proof. So assume
$A_1 \notin E_\ga( 0)$. We would like to
build $A_0 \in E_\ga( 0)$. Set
\begin{align*}
&S=j(T) \big( \ordered{\ga, E( 0)} \big),
\\
&A_0 = \Lev_0(S) \setminus A_1.
\end{align*}
We construct $T^*$:
\begin{align*}
&\Lev_0(T^*)=A_0 \union A_1,
\\
&\forall\gms_1\in A_1 \ T^*_{\ordered{\gms_1}} = T_{\ordered{\gms_1}}.
\end{align*}
We are left with the construction of $T^*_{\ordered{\gms_0}}$ for
$\gms_0 \in A_0$.
For all $\gms \in A_0$ set
\begin{align*}
&A_{\ordered{\gms_0},0}=\Suc_S \big(\ordered{\gms_0} \big),
\\
&A_{\ordered{\gms_0},1}=
        \setof{\ordered{\gms_1, T(\gms_1)_{\ordered{\gms_0}}}}
              {\ordered{\gms_1}\in T, \ \ordered{\gms_0}\in T(\gms_1)},
\\
&\Suc_{T^*} \big( \ordered{\gms_0} \big)=A_{\ordered{\gms_0},0} \union
                A_{\ordered{\gms_0},1},
\\
&\forall\gms_1\in A_{\ordered{\gms_0},1} \ T^*_{\ordered{\gms_0,\gms_1}}=
        T_{\ordered{\gms_1}}.
\end{align*}
We continue one more
level and hope this will convince the reader we indeed can complete $T^*$.
We are left with the construction of $T^*_{\ordered{\gms_0,\gms_1}}$ for
$\ordered{\gms_0,\gms_1}\in A_0 \times A_{\ordered{\gms_0},0}$.
For all $\ordered{\gms_0,\gms_1}\in A_0 \times A_{\ordered{\gms_0},0}$
set
\begin{align*}
&A_{\ordered{\gms_0,\gms_1}, 0}=\Suc_S \big( \ordered{\gms_0,\gms_1} \big),
\\
&A_{\ordered{\gms_0,\gms_1}, 1}=\setof{\ordered{\gms_2,
        T(\gms_2)_{\ordered{\gms_0,\gms_1}}}}
        {\ordered{\gms_2} \in T,\ordered{\gms_0,\gms_1}\in T(\gms_2)},
\\
&\Suc_{T^*} \big( \ordered{\gms_0,\gms_1} \big)=A_{\ordered{\gms_0,\gms_1},0}
\union A_{\ordered{\gms_0,\gms_1},1},
\\
&\forall\gms_2\in A_{\ordered{\gms_0,\gms_1},1}\ 
        T^*_{\ordered{\gms_0,\gms_1,\gms_2}} = T_{\ordered{\gms_2}}.
\end{align*}
We are left with the construction of
$T^*_{\ordered{\gms_0,\gms_1,\gms_2}}$ for
$\ordered{\gms_0,\gms_1,\gms_2}\in A_0 \times
A_{\ordered{\gms_0},0} \times A_{\ordered{\gms_0,\gms_1},0}$ and
we hope that by now the continuation is clear.
\end{itemize}
\end{proof}
\begin{lemma} \label{fill-missing}
Let $\xi_0 < \len(\Es_\ga)$, $T$ a tree such that $\Lev_0(T) \in 
E_\ga( \xi_0)$
and $\forall\ordered{\gms}\in T$
$T_{\ordered{\gms}}$ is an $\Es_\ga$-tree.
Then there is an $\Es_\ga$-tree, $T^*$, satisfying
\begin{enumerate}
\item $\forall \ordered{\gns}\in T \intersect T^*$  
      $T^*_{\ordered{\gns}} \leq T_{\ordered{\gns}}$,
      $T^*(\gns) \leq T(\gns)$,
\item If $p
        \leq
      \bordered{\ordered{\Es_\ga,\ordered{}}, T^*}$ then there is 
        $\ordered{\gms}\in T^* \intersect T$
       such that
        \begin{align*}
       p
        \compatible
      \big(
        \bordered{\ordered{\Es_\ga,\ordered{}}, T^*}
                                \big)_{\ordered{\gms}}.
        \end{align*}
\end{enumerate}
\end{lemma}
\begin{proof}
Our induction hypothesis is that this lemma is
true for $\gms$'s with $\len(\gms) < \len(\Es_\ga)$. The previous lemma is the
case for $\len(\gms)=2$.

Let $S=j(T)(\Es_\ga \restricted \gx_0)$.
The tree $S$ is an $\Es_\ga \restricted \xi_0$-tree. We extend it
step by step
to a full $\Es_\ga$-tree as requested.

Let
\begin{align*}
& A_{\xi_0} = \Lev_0 (T),
\\
& A_{\upto \xi_0} = \Lev_0 (S) \setminus A_{\gx_0}.
\end{align*}
For $\xi_0 < \xi < \len(\Es_\ga)$ do the following:
\begin{align*}
& N_\gx = \Ult(V, E_\ga(\gx)),
\\
& k_\gx([h]_{E_\ga(\gx)}) = j(h)(\Es_\ga \restricted \gx).
\end{align*}
$$
\begin{diagram}
\node{V} \arrow[1]{s} \arrow[1]{e,t}{j} \node[1]{M}
\\
\node{N_\gx=\Ult(V,E_\ga(\gx))} \arrow[1]{ne,b}{{k_\gx}}
\end{diagram}
$$
and recall that $\crit(k_\gx) = (\gk^{++})^{N_\gx}$.
As
\begin{align*}
& \ordof{\ga,E(0),\dotsc,E(\gt),\dotsc} {\gt < \gx} =
k_\gx([id]_{E_\ga(\gx)}) \in \ran(k_\gx)
\end{align*}
there is in $N_\gx$ a preimage for it:
\begin{align*}
& \ordof{\ga',E'(0),\dotsc,E'(\gt') \dotsc} {\gt' < \gx'}.
\end{align*}
As $A_{\gx_0} \in E_\ga(\gx_0)$, $\xi_0 < \gx$ we have 
$\gt' < \gx'$
such that $A_{\gx_0} \in E'_{\ga'}(\gt')$, (where
$\ga'=[\gk(id)]_{E_\ga(\gx)}$. Taking a function $h_\gx$ such that
$[h_\gx]_{E_\ga(\gx)} = \gt'$ we get
\begin{align*}
\setof{\gms} {A_{\gx_0} \intersect \gk^0(\gms) \in
        \gms(h_\gx(\gms))_{\gk(\gms)}}
                \in E_\ga(\gx).
\end{align*}
For each $\gms_0 \in A_{\gx_0}$ we set
\begin{align*}
& A_{\gx, \ordered{\gms_0, T(\gms_0)}} = \setof {\ordered{\gms_1,R} \in
                        T_{\ordered{\gms_0}} }
                    {A_{\gx_0} \intersect \gk^0(\gms_1) \in
                        \gms_1(h_\gx(\gms_1))_{\gk(\gms)}},
\\
& A'_\gx = \sideset{}{^0}\dintersect_{\ordered{\gms_0, T(\gms_0)}}
                A_{\gx, \ordered{\gms_0, T(\gms_0)}}.
\end{align*}
For any tree $R$ which appears in a pair $\ordered{\gms_1, R} \in A'_\gx$
we can invoke by induction our lemma and generate $R^*$ which is a
$\gms_1$-tree. Define now $A_\gx$ as:
\begin{align*}
\ordered{\gms_1, R^*} \in A_\gx  \iff \ordered{\gms_1, R} \in A'_\gx.
\end{align*}
When we have $\setof{A_\gx}{\gx_0 < \gx < \len(\Es_\ga)}$ we set
\begin{align*}
& A_{\downto \gx_0} =  \bigunion_{\gx_0 < \gx < \len(\Es_\ga)} A_\gx
                \setminus (A_{\upto \gx_0} \union A_{\gx_0}),
\\
& \Lev_0(T^*) = A_{\upto \gx_0} \union A_{\gx_0} \union A_{\downto \gx_0},
                        \\
& \forall \ordered{\gms_0} \in A_{\gx_0} \ 
                T^*_{\ordered{\gms_0}} = T_{\ordered{\gms_0}},
\\
& \forall \ordered{\gms_1} \in A_{\downto \gx_0} \ 
         T^*_{\ordered{\gms_1}} = \bigintersect T_{\ordered{\gms_0, \gms_1}}.
\end{align*}
We are left to define $T_{\ordered{\gms_0}}$ for $\ordered{\gms_0} \in 
A_{\upto \gx_0}$.
For each $\gms_0 \in A_{\upto \gx_0}$ set:
\begin{align*}
& A_{\ordered{\gms_0}, \gx_0} =
    \setof {\ordered{\gms_1,R_{\ordered{\gms_0}}}}
           {\ordered{\gms_1,R} \in A_{\gx_0},\ \ordered{\gms_0} \in R},
\\
& A_{\ordered{\gms_0},\upto \xi_0} = \Suc_S (\ordered{\gms_0}) \setminus
        A_{\ordered{\gms_0},\gx_0}.
\end{align*}
For each $\gms_1 \in A_{\gx_0}$ we set
\begin{align*}
& A_{\ordered{\gms_0}, \gx, \ordered{\gms_1,T(\gms_1)_{\ordered{\gms_0}}}} =
        \setof
        {  \ordered{\gms_2, R} \in T_{\ordered{\gms_1}}  }
        { \begin{aligned}
            A_{\gx_0} \intersect \gk^0(\gms_2) \in
                     \gms_2(h_\gx(\gms_2))_{\gk(\gms_2)}   
           \end{aligned}    },
\\
& A'_{\ordered{\gms_0},\gx} =
       \sideset{}{^0}\dintersect_{\ordered{\gms_1,T(\gms_1)_{\ordered{\gms_0}}}}
                A_{\ordered{\gms_0}, \gx, \ordered{\gms_1, T(\gms_1)_{\ordered{\gms_0}}}}.
\end{align*}
We define
\begin{align*}
\ordered{\gms_0, R^*} \in A_{\ordered{\gms_0}, \gx} \iff
\ordered{\gms_0, R} \in A'_{\ordered{\gms_0}, \gx}
\end{align*}
where $R^*$ is generated from R using the current lemma by induction.
Now we set
\begin{align*}
& A_{\ordered{\gms_0},\downto \gx_0} =
              \bigunion_{\gx_0 < \gx < \len(\Es_\ga)} A_{\ordered{\gms_0}, \gx}
        \setminus
                (A_{\ordered{\gms_0}, \upto \gx_0} \union
                 A_{\ordered{\gms_0}, \gx_0}),
\\
& \Suc_{T^*} (\ordered{\gms_0}) = A_{\ordered{\gms_0}, \upto \gx_0} \union
                                A_{\ordered{\gms_0}, \gx_0} \union
                                A_{\ordered{\gms_0}, \downto \gx_0},
\\
& \forall\ordered{\gms_1} \in A_{\gms_0, \gx_0}\ 
        T^*_{\ordered{\gms_0, \gms_1}} = T_{\ordered{\gms_1}},
\\
& \forall\ordered{\gms_2} \in A_{\gms_0, \downto\gx_0}\ 
        T^*_{\ordered{\gms_0, \gms_2}} =
                        \bigintersect T_{\ordered{\gms_1, \gms_2}}.
\end{align*}
This leaves us with the definition of $T_{\ordered{\gms_0, \gms_1}}$ for
$\ordered{\gms_0, \gms_1} \in A_{\upto \gx_0} \times
A_{\ordered{\gms_0},\upto \gx_0}$ which is done exactly as in this step.
\end{proof}
\ref{skeleton} is needed in the proof of \ref{DenseHomogen}. Loosly speaking
it says that if ``something'' happens on all extensions which
are taken from $\dom T$, then that ``something'' happens on all
extensions from $T$.

\ref{skeleton} is proved by induction where \ref{skeleton:2} is the
first case.
\begin{lemma} \label{skeleton:2}
Assume $\len(\Es_\ga)=2$ and let $T$ be $\Es_\ga$-tree. Then there is $T^* \leq T$ such that if
\begin{align*}
        p \leq
      \bordered{\ordered{\Es_\ga,\ordered{}}, T^*}
\end{align*}
then there is $\ordered{\gns_1,\dotsc,\gns_n} \in T$ such that
\begin{align*}
        p  \leq^*
      \big(\bordered{\ordered{\Es_\ga,\ordered{}}, T}
                \big)_{\ordered{\gns_1,\dotsc,\gns_n}}.
\end{align*}
\end{lemma}
\begin{proof}
As is usual in this section the proof is done level by level.
Let us set $T^1 = T$ and it is trivially true that if
\begin{align*}
       p \leq^1
      \bordered{\ordered{\Es_\ga,\ordered{}}, T^1}
\end{align*}
then there is $\ordered{\gns_1} \in T^1=T$ such that
\begin{align*}
       p \leq^*
      \big(\bordered{\ordered{\Es_\ga,\ordered{}}, T^1}
                \big)_{\ordered{\gns_1}}.
\end{align*}
We continue to the second level.
Let us set
\begin{align*}
& A_{\ordered{\gns_1,S^1}} = \Suc_{T^1}(\ordered{\gns_1}),
\\
& A_{\ordered{}}  = \dsintersect_{\ordered{\gns_1,S^1}}
        A_{\ordered{\gns_1,S^1}}
\\
& B_{\ordered{}}  =\setof{\ordered{\gns_2}}
                {\len(\gns_2)=0\text{ or }\Lev_0(T^1)\intersect \gk^0(\gns_2)
                                        \in \gns_2},
\\
& \Lev_0(T^{(0)})=A_{\ordered{}} \intersect B_{\ordered{}},
\\
& T^{(0)}_{\ordered{\gns_2}} = \bigintersect_
                {\substack{
                \ordered{\gns_1} \in T^{(0)}(\gns_2)
                }}
                 T^1_{\ordered{\gns_1,\gns_2}},
\\
& T^2 = T^1 \intersect T^{(0)}.
\end{align*}
Let us assume that
\begin{align*}
        p \leq^2
      \bordered{\ordered{\Es_\ga,\ordered{}}, T^2}.
\end{align*}
There are 2 cases to consider here:
\begin{enumerate}
\item $p \leq^2
      \big(\bordered{\ordered{\Es_\ga,\ordered{}}, T^2}
        \big)_{\ordered{\gns_1,\gns_2}}$ where
        $\ordered{\gns_1, \gns_2} \in T^2$:
        At once we have \linebreak[0] $\ordered{\gns_1,\gns_2}\in T^1 \leq T$.
\item $p \leq^2
      \big(\bordered{\ordered{\Es_\ga,\ordered{}}, T^2}
        \big)_{\ordered{\gns_2,\gns_1}}$ where
$\ordered{\gns_2} \in T^2$, $\ordered{\gns_1}\in T^2(\gns_2)$:
        By construction $\forall \ordered{\gms}\in T^2(\gns_2)\ 
                \ordered{\gns_2}\in 
                T^1_{\ordered{\gms}}$. As $\ordered{\gns_1} \in 
                T^2(\gns_2)$ we get $\ordered{\gns_1,\gns_2} \in T^1 \leq T$.
\end{enumerate}
We show how to continue to the third level.
Let us set
\begin{align*}
& A_{\ordered{\gns_1,S^1,\gns_2,S^2}} = \Suc_{T^2}(\ordered{\gns_1,\gns_2}),
\\
& A_{\ordered{\gns_1,S^1}}  = \dsintersect_{\ordered{\gns_2,S^2}}
        A_{\ordered{\gns_1,S^1,\gns_2,S^2}},
\\
& B_{\ordered{\gns_1,S^1}}  = \setof {\ordered{\gns_3}}
                {\len(\gns_3)=0\text{ or }
                 \Suc_{T^2}(\ordered{\gns_1,S^1})\intersect 
                        \gk^0(\gns_3) \in \gns_3},
\\
& \Lev_0(T^{(0)}) = \Lev_0(T^2),
\\
& \Suc_{T^{(0)}}(\ordered{\gns_1})= A_{\ordered{\gns_1,S^1}} \intersect
                                        B_{\ordered{\gns_1,S^1}},
\\
& T^{(0)}_{\ordered{\gns_1,\gns_3}} = \bigintersect_
                        {\ordered{\gns_2} \in 
                                T^{(0)}_{\ordered{\gns_1}}(\gns_3)}
                        T^2_{\ordered{\gns_1,\gns_2,\gns_3}},
\\
& A_{\ordered{}}=\dsintersect_{\ordered{\gns_1,S^1}}
        A_{\ordered{\gns_1,S^1}},
\\
& B_{\ordered{}}=\setof{\ordered{\gns_3}}
                     {\len(\gns_3)=0\text{ or }
                        \Lev_0(T^2)\intersect \gk^0(\gns_3) \in \gns_3},
\\
& \Lev_0(T^{(1)}) = A_{\ordered{}} \intersect B_{\ordered{}},
\\
& \Suc_{T^{(1)}}(\ordered{\gns_3})= \bigintersect_
                        {\ordered{\gns_1, \gns_2} \in 
                        T^{(0)}(\gns_3)}
                        T^2_{\ordered{\gns_1,\gns_2,\gns_3}},
\\
& T^3 = T^2 \intersect T^{(0)} \intersect T^{(1)}.
\end{align*}
Let us assume that
\begin{align*}
        p \leq^3
      \bordered{\ordered{\Es_\ga,\ordered{}}, T^3}.
\end{align*}
There are 3 cases to consider here:
\begin{enumerate}
\item $p \leq^*
      \big(\bordered{\ordered{\Es_\ga,\ordered{}}, T^3}
                \big)_{\ordered{\gns_1,\gns_2,\gns_3}}$
        where
        $\ordered{\gns_1,\gns_2,\gns_3} \in T^3$:
        At once we get \linebreak
                $\ordered{\gns_1,\gns_2,\gns_3}\in \linebreak[0] T$.

\item $p \leq^*
      \big(\bordered{\ordered{\Es_\ga,\ordered{}}, T^3}
                \big)_{\ordered{\gns_1,\gns_3,\gns_2}}$
        where
        $\ordered{\gns_1,\gns_3} \in T^3$,
        $\ordered{\gns_2}\in T^3_{\ordered{\gns_1}}(\gns_3)$:
        In this case
        $\forall \ordered{\gms}\in T^3_{\ordered{\gns_1}}(\gns_3)\ 
                \ordered{\gns_3}\in 
                T^3_{\ordered{\gns_1,\gms}}$. 
        As $\ordered{\gns_2} \in 
                T^3_{\ordered{\gns_1}}(\gns_3)$ we get 
                $\ordered{\gns_3} \in T_{\ordered{\gns_1,\gns_2}}$.

\item $p \leq^*
      \big(\bordered{\ordered{\Es_\ga,\ordered{}}, T^3}
                \big)_{\ordered{\gns_3,\gns_1,\gns_2}}$
        where
        $\ordered{\gns_3} \in T^3$,
        $\ordered{\gns_1, \gns_2} \in T^3(\gns_3)$:
        \linebreak[0]
        $\forall \ordered{\gms_1,\gms_2}\in T^3(\gns_3)\ 
                \ordered{\gns_3}\in 
                T^3_{\ordered{\gms_1,\gms_2}}$
        hence $\ordered{\gns_3}\in 
                T_{\ordered{\gns_1,\gns_2}}$.
\end{enumerate}
In this way we continue to all levels.
\end{proof}
\begin{lemma} \label{skeleton}
Let $T$ be $\Es_\ga$ tree. Then there is $T^* \leq T$ such that if
\begin{align*}
        p \leq
      \bordered{\ordered{\Es_\ga,\ordered{}}, T^*}
\end{align*}
then there is $\ordered{\gns_1,\dotsc,\gns_n} \in T$ such that
\begin{align*}
        p  \leq^*
      \big(\bordered{\ordered{\Es_\ga,\ordered{}}, T}
                \big)_{\ordered{\gns_1,\dotsc,\gns_n}}.
\end{align*}
\end{lemma}
\begin{proof}
The proof is by induction on $\len(\Es_\ga)$. The first case was
done in \ref{skeleton:2}. The proof is almost the same. We just make
sure to invoke the induction hypothesis while repeating the construction.

Construction of $T^1$ and $T^2$ is exactly like in \ref{skeleton:2}.
We show the construction at the 3rd level.

Let us set
\begin{align*}
& A_{\ordered{\gns_1,S^1,\gns_2,S^2}} = \Suc_{T^2}(\ordered{\gns_1,\gns_2}),
\\
& A_{\ordered{\gns_1,S^1}}  = \dsintersect_{\ordered{\gns_2,S^2}}
        A_{\ordered{\gns_1,S^1,\gns_2,S^2}},
\\
& B_{\ordered{\gns_1,S^1}}  = \setof {\ordered{\gns_3}}
                {\len(\gns_3)=0\text{ or }
                 \Suc_{T^2}(\ordered{\gns_1,S^1})\intersect 
                        \gk^0(\gns_3) \in \gns_3},
\\
& \Lev_0(T^{(0)}) = \Lev_0(T^2),
\\
& \Suc_{T^{(0)}}(\ordered{\gns_1})= A_{\ordered{\gns_1,S^1}} \intersect
                                        B_{\ordered{\gns_1,S^1}},
\\
& T^{(0)}_{\ordered{\gns_1,\gns_3}} = \bigintersect_
                        {\ordered{\gns_2} \in 
                                T^{(0)}_{\ordered{\gns_1}}(\gns_3)}
                        T^2_{\ordered{\gns_1,\gns_2,\gns_3}},
\\
& A'_{\ordered{}}=\dsintersect_{\ordered{\gns_1,S^1}}
        A_{\ordered{\gns_1,S^1}},
\\
& A_{\ordered{}}=\setof{\ordered{\gns_3,S^3}}{\ordered{\gns_3,S^{3\prime}}\in 
                        A'_{\ordered{}}\ S\text{ is generated
                                from }S^{3\prime}\text{ by induction}},
\\
& B_{\ordered{}}=\setof{\ordered{\gns_3}}
                     {l(\gns_3)=0\text{ or }
                        \Lev_0(T^2)\intersect \gk^0(\gns_3) \in \gns_3},
\\
& \Lev_0(T^{(1)}) = A_{\ordered{}} \intersect B_{\ordered{}},
\\
& \Suc_{T^{(1)}}(\ordered{\gns_3})= \bigintersect_
                        {\ordered{\gns_1, \gns_2} \in 
                        T^{(0)}(\gns_3)}
                        T^2_{\ordered{\gns_1,\gns_2,\gns_3}},
\\
& T^3 = T^2 \intersect T^{(0)} \intersect T^{(1)}.
\end{align*}
Let us assume that
\begin{align*}
        p_2 \append p_1 \append p_0 = p \leq^3
      \bordered{\ordered{\Es_\ga,\ordered{}}, T^3}.
\end{align*}
There are 3 cases to consider here:
\begin{enumerate}
\item $p \leq^*
      \big(\bordered{\ordered{\Es_\ga,\ordered{}}, T^3}
        \big)_{\ordered{\gns_1,\gns_2,\gns_3}}$
        where
        $\ordered{\gns_1,\gns_2,\gns_3} \in T^3$:
        At once we get  $\linebreak[0]\ordered{\gns_1,\gns_2,\gns_3}\in T$.
\item $p \leq^*
      \big(\bordered{\ordered{\Es_\ga,\ordered{}}, T^3}
        \big)_{\ordered{\gns_1,\gns_3,\gns_2}}$
        where
         $\ordered{\gns_1,\gns_3} \in T^3$,
        $\ordered{\gns_2}\in T^3_{\ordered{\gns_1}}(\gns_3)$:
        In this case
        $\forall \ordered{\gms}\in T^3_{\ordered{\gns_1}}(\gns_3)\ 
                \ordered{\gns_3}\in 
                T^3_{\ordered{\gns_1,\gms}}$. 
        As $\ordered{\gns_2} \in 
                T^3_{\ordered{\gns_1}}(\gns_3)$ we get 
                $\ordered{\gns_3} \in T_{\ordered{\gns_1,\gns_2}}$.
\item $p \leq^2
      \big(\bordered{\ordered{\Es_\ga,\ordered{}}, T^3}
        \big)_{\ordered{\gns_3}}$
        where
        $\ordered{\gns_3} \in T^3$ and $p_2 \append p_1 \leq^2 
              \big(\bordered{\ordered{\gns_3,\ordered{}}, \linebreak[0]
                T^3(\gns_3)}$
        By induction there is $\ordered{\gns_1,\gns_2} \in T(\gns_3)$
        such that 
        \begin{align*}
        \linebreak[0] p_2 \append p_1 \linebreak[0] \leq^*  \linebreak[0]
              \big(\bordered{\ordered{\gns_3,\ordered{}}, \linebreak[0]
                        T^2(\gns_3)}
                        \big)_{\ordered{\gns_1,\gns_2}}.
        \end{align*}
        By construction
        $\forall \ordered{\gms_1,\gms_2}\in T^2(\gns_3)\ 
                \ordered{\gns_3}\in 
                T^2_{\ordered{\gms_1,\gms_2}}$
        hence $\ordered{\gns_3}\in 
                T_{\ordered{\gns_1,\gns_2}}$.
\end{enumerate}
In this way we continue to all levels.
\end{proof}
%
%       F O R C I N G   - Single chunk
%
\section{$P_{\Es}$-Forcing} \label{PEForcing}
\begin{definition}
A condition in $P_{\Es}^*$ is of the form
\begin{equation*}
\setof{\ordered{\ggs, p^\ggs}}{\ggs \in g} \union
        \set{T}
\end{equation*}
where
\begin{enumerate}
\item $g \subseteq \Es$, $\power{g} \leq \kappa$,
\item $\min \Es\in g$ and $g$ has a maximal element,
\item $p^\ggs \in V_\gk$ is an extender sequence. We allow $p^\ggs=\emptyset$,
\item $p^0 = (p^{\max g})^0$. \par\noindent
        This condition is not really needed here. It 
        is needed in a later forcing based on this one,
\item $T$ is a $\max g$-tree such that for all $t \in T$
        $p^{\max g}\append t$ is $^0$-increasing,
\item For all $\ggs \in g$, $p^{\max g}$ is not permitted
                to $p^\ggs$,
\item $\forall \ordered{\gns}\in T \ 
        \power{ \setof{\ggs \in g}{\gns \text{ is permitted to } p^\ggs }}
                                                \leq \gk^0(\gns)$,
\item $\forall \ordered{\gns}\in T$ if $\gns$ is permitted to
        $p^\gbs,p^\ggs$ then 
        $\gp_{\max g,\gbs}(\gns) \not =
         \gp_{\max g, \ggs}(\gns)$.
\end{enumerate}
We write $\mc(p)$, $p^{\mc}$, $T^p$, $\Es(p)$, $\supp p$ for
$\max g$, $p^{\max g}$, $T$, $\Es$, $g$ respectively.
\end{definition}
\begin{figure}[htb]
$$
\begin{diagram}
\\
\node[4]{} \node{\gms_0} \node[5]{\gms_1} \node[5]{\gms_2} \node[5]{\gms_3} \node[5]{\gms_4} \node[5]{T} 
\\
\node[4]{} \arrow[22]{e,-} \node[5]{}
\\
\node{\text{Support}} \node[3]{} \node{\Es_\gk} \node[5]{\Es_{\ga_1}} \node[5]{\Es_{\ga_2}} \node[5]{\Es_{\ga_3}} \node[5]{\Es_{\ga_4} = \mc} \node[5]{} 
\end{diagram}
$$
\caption{An Example of Condition in $P^*_{\Es}$}
\end{figure}
%\fi
%
%  Order on single chunk
%
\begin{definition}
Let $p, q \in P^*_{\bar{E}}$. We say that $p$ is a Prikry extension of $q$
($p \leq^* q$ or $p \leq^0 q$) if
\begin{enumerate}
\item $\supp p \supseteq \supp q$,
\item $\forall \gga \in \supp q \ p^\gga=q^\gga$,
\item $T^p \leq T^q$.
\end{enumerate}
\end{definition}
We include in this definition the degenerate case $\len(\Es)=0$. There is
neither extender nor tree in this case. By $p \leq^* p$  we
mean $p = q$.
\begin{figure}[htb]
$$
\begin{diagram}
\node[4]{} \node{\gms_0} \node[2]{\gns_0} \node[3]{\gms_1} \node[2]{\gns_1} \node[3]{\gms_2} \node[5]{\gms_3} \node[5]{\gms_4} \node[5]{\gns_2} \node[5]{\gp_{\gb_2,\ga_4}^{-1}T}
\\
\node[4]{} \arrow[27]{e,-} \node[5]{}
\\
\node{\text{Support}} \node[3]{} \node{\Es_{\gk}} \node[2]{\Es_{\gb_0}} \node[3]{\Es_{\ga_1}} \node[2]{\Es_{\gb_1}} \node[3]{\Es_{\ga_2}} \node[5]{\Es_{\ga_3}} \node[5]{\Es_{\ga_4}} \node[5]{\Es_{\gb_2}=\mc} 
\end{diagram}
$$
\caption{An Example of Direct Extension}
\end{figure}
%
% The actual forcing
%
\begin{definition}
A condition in $P_{\Es}$ is of the form
\begin{equation*}
p_n \append \dotsb \append p_0
\end{equation*}
where
\begin{itemize}
\item $p_0 \in P^*_{\Es}$,
\item $p_1 \in P^*_{\gms_1}$,
\item $\vdots$,
\item $p_n \in P^*_{\gms_n}$,
\end{itemize}
where $\Es, \gms_1, \dotsc, \gms_n$ are extender sequence
systems satisfying
\begin{equation*}
\gk^0(\gms_n) < \dotsb < \gk^0(\gms_1) < \gk^0(\bar{E}).
\end{equation*}
\end{definition}
\begin{figure}[htb]
$$
\begin{diagram}
\node{\gt_0} \node[2]{\gt_1} \node[2]{\gt_2} \node[2]{\gt_3} \node[3]{\gt_4} \node[2]{R} \node[1]{}
\node{\gms_5} \node[2]{\gms_6} \node[2]{\gms_7} \node[2]{\gms_8} \node[3]{\gms_9} \node[2]{S} \node[1]{}
\node{\gns_0} \node[2]{\gns_1} \node[2]{\gns_5} \node[2]{\gns_6} \node[3]{\gns_4} \node[2]{T} 
\\
\arrow[11]{e,-} \node[13]{}
\arrow[11]{e,-} \node[13]{}
\arrow[11]{e,-} \node[2]{}
\\
\node{\gms_0} \node[2]{\gms_1} \node[2]{\gms_2} \node[2]{\gms_3} \node[3]{\gms_4=\mc} \node[2]{} \node[1]{}
\node{\gns_0} \node[2]{\gns_1} \node[2]{\gns_2} \node[2]{\gns_3} \node[3]{\gns_4=\mc} \node[2]{} \node[1]{}
\node{\Es_{\gk}} \node[2]{\Es_{\ga_1}} \node[2]{\Es_{\ga_2}} \node[2]{\Es_{\ga_3}} \node[3]{\Es_{\ga_4}=\mc} \node[2]{} 
\end{diagram}
$$
\caption{An Example of a Condition in $P_{\Es}$}
\end{figure}
\begin{definition}
Let $p,q \in P_{\Es}$. We say that $p$ is a Prikry extension of $q$
($p \leq^* q$ or $p \leq^0 q$)
if $p,q$ are of the form
\begin{equation*}
\begin{split}
p &= p_n \append \dotsb \append p_0,
\\
q &= q_n \append \dotsb \append q_0,
\end{split}
\end{equation*}
and
\begin{itemize}
\item $p_0,q_0 \in P^*_{\Es},\ p_0 \leq^* q_0$,
\item $p_1,q_1 \in P^*_{\gms_1},\ p_1 \leq^* q_1$,
\item $\vdots$,
\item $p_n,q_n \in P^*_{\gms_n},\ p_n \leq^* q_n$.
\end{itemize}
\end{definition}
%
% 1-point extension
%
\begin{definition} \label{dfn:enlarge}
Let $p \in P^*_{\bar{E}}$ and $\ordered{\gns} \in T^p$. We =define
$(p)_{\ordered{\gns}}$ to be
$p'_1 \append p'_0$ where
\begin{enumerate}
\item $\supp p'_1=
                \setof{\pi_{\mc(p),\ggs}(\gns)}
              {\ggs \in \supp p, \, \gk^0(p^\ggs) < \gk^0(\gns)}$,
\item $p_1^{\prime\pi_{\mc(p),\ggs}(\gns)}=p^\ggs$,
\item $T^{p'_1} = T^p(\gns)$,
\item $\supp p'_0 = \supp p$,
\item $\forall \ggs \in \supp p'_0 \ p_0^{\prime\ggs}=
        \begin{cases}
        \pi_{\mc(p),\ggs}(\gns) & 
                                \gk^0(p^\ggs) < \gk^0(\gns)
        \\
        p^\ggs  &       \text{otherwise}
        \end{cases}$,
\item $T^{p'_0}=T^p_{\ordered{\gns}}$.
\end{enumerate}
\end{definition}
\begin{definition}
Let $p,q \in P_{\Es}$. We say that $p$ is a $1$-point extension of $q$
($p \leq^1 q$) if
$p,q$ are of the form
\begin{equation*}
\begin{split}
p &= p_{n+1} \append p_n \append \dotsb \append p_0,
\\
q &= q_n \append \dotsb \append q_0,
\end{split}
\end{equation*}
and there is $0 \leq k \leq n$ such that
\begin{itemize}
\item $p_i,q_i \in P^*_{\gms_i},\ p_i \leq^* q_i$ for $i=0,\dotsc,k-1$,
\item $p_{i+1},q_i \in P^*_{\gms_i},\ p_{i+1} \leq^* q_i$ for $i=k+1,\dotsc,n$,
\item There is $\ordered{\gns} \in T^{q_k}$ such that
        $p_{k+1} \append p_k \leq^* (q_k)_{\ordered{\gns}}$.
\end{itemize}
\end{definition}
\begin{figure}[htb]
$$
\begin{diagram}
\node{\gms_0} \node[4]{\gms_1} \node[6]{\gms_3} \node[4]{\gms_4} \node[2]{T(\gns)} \node[2]{}
\node{\gns^0} \node[4]{\gp_{\ga_4,\ga_1}(\gns)} \node[4]{\gms_2} 
			\node[4]{\gp_{\ga_4,\ga_3}(\gns)} \node[4]{\gns}
		\node[2]{T_{\ordered{\gns}}}
\\
\arrow[16]{e,-} \node[18]{}
\arrow[19]{e,-} \node[18]{}
\\
\node{\gns^0} \node[4]{\gp_{\ga_4,\ga_1}(\gns)} \node[6]{\gp_{\ga_4,\ga_3}(\gns)} \node[4]{\gns} \node[4]{}
\node{\Es_{\gk}} \node[4]{\Es_{\ga_1}} \node[4]{\Es_{\ga_2}} \node[4]{\Es_{\ga_3}} \node[4]{\Es_{\ga_4}}
\end{diagram}
$$
\caption{An Example of 1-point extension}
\end{figure}
%
% n-point extension
%
\begin{definition}
Let $p,q \in P_{\Es}$. We say that $p$ is an $n$-point extension of $q$
($p \leq^n q$) if there are $p^n, \dotsc, p^0$ such that
\begin{align*}
        p=p^n \leq^1 \dotsb \leq^1 p^0=q.
\end{align*}
\end{definition}
\begin{definition}
Let $p,q \in P_{\Es}$. We say that $p$ is an extension of $q$
($p \leq q$) if there is n such that $p \leq^n q$.
\end{definition}
\par\noindent Later on by $\PE$ we mean $\ordered{\PE,\leq}$.
\begin{note}
When $\len(\Es)=1$ the forcing $P_{\Es}$ is
the Gitik-Magidor forcing from section 1 of \cite{Moti}.
When $\len(\Es) < \gk$ the forcing $P_{\Es}$ is similar to the forcing
defined in \cite{Miri}.
\end{note}
In several places we want to prevent enlargment of the support of
a condition. This makes all the conditions which are stronger than
some condition but with the same support resemble Radin forcing.
The following definition catches the meaning of not enlarging the
support. The `resemblence' we look for is \ref{GenericForRadin}.
\begin{definition}
Let $p,q \in P_{\Es}$. We say that $p \leq^*_R q$ 
if
\begin{enumerate}
\item $p \leq^* q$,
\item $\supp p = \supp q$.
\end{enumerate}
\end{definition}
\begin{definition}
Let $p,q \in P_{\Es}$. We say that $p \leq^1_R q$ 
if
\begin{enumerate}
\item $p \leq^1 q$,
\item In the definition of $\leq^1$ we can substitute $\leq^*$ by
        $\leq^*_R$.
\end{enumerate}
\end{definition}
\begin{definition}
Let $p,q \in P_{\Es}$. We say that $p \leq^n_R q$ if 
there are $p^n, \dotsc, p^0$ such that
\begin{align*}
        p=p^n \leq^1_R \dotsb \leq^1_R p^0=q.
\end{align*}
\end{definition}
\begin{definition}
Let $p,q \in P_{\Es}$. We say that $p \leq_R q$ if there is $n$
such that $p \leq^n_R q$.
\end{definition}
\begin{note}
The above definitions imply that
if $q \leq p$ then there is $r$ such that $q \leq^* r \leq_R p$.
\end{note}
\begin{definition}
Let $\ges$ be an extender sequence such that $\gk^0(\ges)< \gk^0(\Es)$.
\begin{align*}
\PE/\Pe = \setof {p}{q\in \Pe, q\append p \in \PE}.
\end{align*}
\end{definition}
%
%
% Basic Properties
%
\section{Basic Properties of $P_{\Es}$} \label{BasicProperties}
\begin{claim}
$P_\Es$ satisfies $\gk^{++}$-c.c.
\end{claim}
\begin{proof}
The usual $\Delta$-lemma argument on the support will do.
\end{proof}
\begin{claim} \label{SubForcing}
Let $p \in \PE$, $P^* = \setof {q \leq_R p} {p \in \PE}$. Then 
\begin{enumerate}
\item $\ordered{P^*,\leq_R}$ satisfies $\gk^{+}$-c.c.,
\item $\ordered{P^*,\leq_R}$ is sub-forcing of $\ordered{\PE/p,\leq}$.
\end{enumerate}
\end{claim}
\begin{proof}
Showing $\gk^{+}$-c.c. is trivial.

Showing that $P^*$ is sub-forcing of $\PE/p$ amounts to showing that 
any maximal anti-chain of
$P^*$ is also a maximal anti-chain of $\PE/p$.

Let $A$ be a maximal anti-chain of $P^*$.
Let $q \in \PE/p$. As $q \leq p$, there is
$r' \in P^*$ such that $q \leq^* r' \leq_R p$.
Assume that $r' = r'_n \append \dotsb \append r'_0$. Then also
$q = q_n \append \dotsb \append q_0$. Let $r_i$ be $r'_i$ with 
$T^{r'_i}$ substituted by $T^{r'_i} \intersect 
\gp_{\mc(q_i),\mc(r'_i)}(T^{q_i})$ and $r = r_n \append \dotsb \append r_0$.
As $r \in P^*$ and $A$ is a maximal anti-chain there
is $a \in A$ such that $a \compatible r$. Take $s \leq_R a, r$.
Considering how we constructed $r$ from $r'$ we must have
$t \leq^* s$ such that $t \leq q$. Hence $q \compatible a$.
So we get that $A$ is a maximal anti-chain of $\PE/p$.
\end{proof}
\begin{claim} \label{GenericForRadin}
Let $p \in \PE$, $P^* = \setof {q \leq_R p} {p \in \PE}$. 
Then there is $r \in R_{\mc(p)}$ such that $P^* \simeq R_{\mc(p)}/r$
\end{claim}
\begin{proof}
For simplicity assume that $p = p_0$. Then we set
$r = \ordered{\ordered{\mc(p_0), p_0^{\mc}}, \linebreak[0] T^p}$.

We give the isomorphism: The image of $q \in P^*$ is
$s \in \ordered{R_{\mc(p)}/r, \leq}$ such that
\begin{enumerate}
\item $q \leq^* s$,
\item $T^{s_i} = T^{q_i}$ where
        $s = s_n \append \dotsb \append s_0$,
        $q = q_n \append \dotsb \append q_0$.
\end{enumerate}
\end{proof}
\par\noindent
Let $G$ be $\PE$-generic.
\par\noindent
\begin{definition}
$\Es_G$ is the enumeration of
        $\setof {\Es(p_k)} {p_n \append \dotsb \append p_0 \in G}$
 ordered increasingly by $\gk^0(\Es(p_k))$.
\end{definition}
\begin{definition} Let $\gz < \otp(\Es_G)$. Then
\begin{enumerate}
\item $G \restricted \gz = \setof {p_n \append \dotsb \append p_k}
                            {p_n \append \dotsb \append p_k \append
                                        \dotsb \append p_0 \in G, 
                                \Es(p_k) = \Es_G(\gz)}$,
\item $G \setminus \gz = \setof {p_{k-1} \append \dotsb \append p_0}
                            {p_n \append \dotsb \append p_k \append \dotsb 
                                        \append p_0 \in G, \linebreak[0]
                                \Es(p_k) = \Es_G(\gz)}$.
\end{enumerate}
\end{definition}
\begin{definition}
\begin{align*}
& M^{\gas}_G = 
        \begin{cases}
        \bigunion \setof{M^{\gms}_G}{p \in G, \, \gms = p^{\gas}} \union 
                                    \set {\gas} &
                        \exists p\in G\ \gas\in \supp p
        \\
        \set{\gas}  & \text{otherwise}
        \end{cases}
\\
& C^{\gas}_G = \setof
        {\gk(\gms)}
        {\gms \in M^{\gas}_G}
\end{align*}
\end{definition}
\begin{proposition}
\begin{enumerate}
\item $C^{\Es_\gk}_G\setminus\set{\gk}$ is a club in $\gk$,
\item $C^{\Es_\ga}_G\setminus\set{\ga}$ is unbounded in $\gk$,
\item $\gas \not= \gbs \implies C^{\gas}_G \not= C^{\gas}_G$.
\end{enumerate}
\end{proposition}
\begin{proof}
The first two claims are immediate as these are sequences which
are generated by Radin forcing.

The last is by density and noticing that when $p^{\gas},p^{\gbs}$ are
permitted for $\gns$ we required $\gp_{\mc(p),\gas}(\gns) \not=
\gp_{\mc(p),\gbs}(\gns)$.
\end{proof}
%
%
% Homogeneity of Dense Open sets
%
\section{Homogeneity in Dense Open Subsets} \label{HomogenDense}
Our aim in this section is to prove the following
\begin{theorem} \label{DenseHomogen}
Let $D \subseteq \PE$ be dense open and $p=p_k\append\dotsb\append p_0 \in \PE$. 
Then there is $p^* \leq^* p$
such that 
\begin{multline*}
\exists S^k\ \exists n_k \forall \ordered{\gns_{k,1},\dotsc,\gns_{k,{n_k}}} 
                \in S^k \, \dotsc
\exists S^0\ \exists n_0 \forall \ordered{\gns_{0,1},\dotsc,\gns_{0,{n_0}}}
                \in S^0 \,
\\
         (p^*_k)_{\ordered{\gns_{k,1},\dotsc,\gns_{k,n_k}}} \append
          \dotsb \append
         (p^*_0)_{\ordered{\gns_{0,1},\dotsc,\gns_{0,n_0}}}
                         \in D
\end{multline*}
where
\begin{enumerate}
\item $S^i \subseteq T^{p^*_i} \restricted [V_\gk]^{n_i}$,
\item $\forall l<n_i \,\forall \ordered{\gns_1,\dotsc,\gns_{l}} \in S^i \,
        \exists \gx \ 
        \Suc_{S^i}(\ordered{\gns_1,\dotsc,\gns_{l}}) 
                                                \in E_{\mc(p^*_i)}(\gx)$.
\end{enumerate}
\end{theorem}
\par\noindent
The proof is done by a series of lemmas.
\begin{definition} \label{EnlargeEvery}
Let $p \in P^*_{\Es}$.
Let $s$ be a function such that $\dom s \subseteq \Es$.
and for all $\gas,\gbs \in \dom s$, $\gas \neq \gbs$
\begin{enumerate}
\item $s(\gas)$ is an extender sequence,
\item $\len(s(\gas))=\len(s(\gbs))$,
\item $\gk^0(s(\gas)) = \gk^0(s(\gbs))$,
\item $s(\gas) \not= s(\gbs)$.
\end{enumerate}
We define
$(p)_{\ordered{s}}$ to be $p'_1 \append p'_0$ where
\begin{enumerate}
\item $\supp p'_1=
                \setof{s(\gas)}
              {\gas \in \supp p\intersect\dom s, \, \gk^0(p^\gas) < \gk^0(s(\gas))}$,
\item $\forall \gas\in \supp p'_1\ p_1^{\prime s(\gas)}=
                p^\gas$,
\item If $s(\mc(p))\in T^{p}$ then $T^{p'_1}=T^{p}(s(\mc(p)))$. Otherwise
      we leave $T^{p'_1}$ undefined,
\item $\supp p'_0 = \supp p$,
\item $\forall \gas \in \supp p'_0 \ p_0^{\prime\gas}=
        \begin{cases}
        s(\gas) &
                                \gas\in\dom s\text{ and }
                                \gk^0(p^\gas) < \gk^0(s(\gas))
        \\
        p^\gas  &       \text{otherwise}
        \end{cases}$,
\item If $s(\mc(p))\in T^{p}$ then $T^{p'_0}=T^{p}_{\ordered{s(\mc(p))}}$.
        Otherwise we leave $T^{p'_0}$ undefined.
\end{enumerate}
\end{definition}
\begin{definition}
Let $p \in P^*_{\Es}$.
Let $s$ be a function with $\dom s = {1,\dotsc,n}$ such that for all $i$
$s(i)$ satisfies definition \ref{EnlargeEvery}. Then we define
$(p)_{\ordered{s}}$ as $p^n$ where $p^n$ is defined by induction as
follows:
\begin{align*}
& p^0 = p,
\\
& p^{i+1} = p^i_i \append \dotsb \append p^i_1 \append 
        (p^i_0)_{\ordered{s(i+1)}}.
\end{align*}
\end{definition}
\par\noindent
We note the following: If $\ordered{\gns_1,\dotsc,\gns_n} \in T^p$ and we
set for all $1\leq i \leq n$
\begin{align*}
s(i)=\setof{\ordered{\gas, \gp_{\mc(p),\gas}(\gns_i)}}{\gas \in \supp p}
\end{align*}
then
\begin{align*}
(p)_{\ordered{\gns_1,\dotsc,\gns_n}}=(p)_{\ordered{s}}.
\end{align*}

We use this operation also in cases where $p$ is not strictly a
condition. That is if $p \union \set{T} \in \PE$ we also
use $(p)_{\ordered{s}}$. In this case we ignore the trees in
the definition.

This definition is used in the proof of the homogeneity for
the following reason:
Beforehand we do not know what a legitimate extension is. By checking
with all the possible $\gm$'s we check on all legitimate conditions
which {\em might} be extensions.
%
% Canonical  for Dense Open - 1nd level
%
\begin{claim} \label{canon-dense:n}
Let $D$ be dense open in $\PE/\Pe$, $p = p_0 \in \PE/\Pe$, $0<n<\gw$.
Then there is $p^* \leq^* p$
such that one and only one of the following is true:
\begin{enumerate}
\item There is $S \subseteq T^{p^*}\restricted [V_\gk]^{n}$ such that 
\begin{enumerate}
\item $\forall k<n\,\exists \gx<\len(\Es)\ 
                \Suc_S(\ordered{\gns_1,\dotsc,\gns_k})
                         \in E_{\mc(p^*)}(\gx)$,
\item $\forall \ordered{\gns_1,\dotsc,\gns_{n}} \in S$
                $(p^*)_{\ordered{\gns_1,\dotsc,\gns_{n}}} \in D$.
\end{enumerate}
\item $\forall \ordered{\gns_1,\dotsc,\gns_{n}} \in 
                        T^{p^*}
        \forall q \leq^* (p^*)_{\ordered{\gns_1,\dotsc,\gns_{n}}}\ 
               q \notin D$.
\end{enumerate}
\end{claim}
\begin{proof}
We give the proof for $n=1$. Adapting the proof for higher $n$'s require
that whenever we enumerate singletons we should enumerate $n$-tuples
and when we use $j$ we should use $j_n$.

We start an induction on $\gx$ in which we build
\begin{align*}
\ordof{\gas^\gx, u^\gx}{\gx < \gk}.
\end{align*}
We start by setting
\begin{align*}
u^0 &= p_0 \setminus \set{T^{p_0}},
\\
\gas^0 &= \mc(p_0),
\\
T^0 &= T^{p_0} \restricted \pi^{-1}_{\gas^0,0}\setof{\gns}
                                     {\gk^0(\gns) \text{ is inaccessible}},
\end{align*}
and taking an increasing enumeration
\begin{align*}
\setof{\gk^0(\gns)}{\ordered{\gns} \in T^0}
                = \ordof{\gt_\gx}{\gx < \gk}.
\end{align*}
Assume that we have constructed
\begin{align*}
\ordof{\gas^\gx, u^\gx}{\gx < \gx_0}.
\end{align*}
\par\noindent
We have $2$ cases:
$\gx_0$ is limit: Choose $\gas^{\gx_0} > \gas^\gx$ for all $\gx < \gx_0$
and set
\begin{align*}
u^{\gx_0} &= \bigunion_{\gx < \gx_0} u^\gx \union
        \set{\ordered{\gas^{\gx_0}, t}}
                \text{ where }\gk^0(t) = \gt_{\gx_0}.
\end{align*}
\par\noindent
$\gx_0 = \gx + 1:$
For each $\gns$ such that
$\gk^0(\gns)=\gt_\gx$ we set
\begin{align*}
 S(\gns) = 
        \big( \prod_{\substack{\gas \in \supp u^\gx\\
                        \gk((u^\gx)^{\gas})<\gt_\gx}} 
                \setof {\gms}{\gk^0(\gms)=\gk^0(\gns)} \big)
                                     \times
         \set {\ordered{\gns}}.
\end{align*}
Let
\begin{align*}
S = \bigunion_{\gk^0(\gns)=\gt_\gx} 
        S(\gns)
\end{align*}
and set enumeration of $S$
\begin{align*}
S=\ordof{s^{\gx_0,\gr}} {\gr < \gt_{\gx_0}}.
\end{align*}
There are fewer than $\gt_{\gx_0}$ elements in $S$. We use $\gt_{\gx_0}$ 
as this is the maximum size $S$ can have which is not `killing' the induction.

We do induction on $\gr$ which  builds
\begin{align*}
\ordof{\gas^{\gx_0,\gr}, u_0^{\gx_0,\gr},
                        T_0^{\gx_0,\gr},u_1^{\gx_0,\gr},T_1^{\gx_0,\gr}}
                                                {\gr < \gt_{\gx_0}}
\end{align*}
from which we build $\ordered{\gas^{\gx_0}, u^{\gx_0}}$.
Set
\begin{align*}
& \gas^{\gx_0,0}=\gas^\gx,
\\
& u^{\gx_0,0}_0 = u^\gx_0.
\end{align*}
Assume we have constructed
$\ordof{\gas^{\gx_0,\gr}, u_0^{\gx_0,\gr},T_0^{\gx_0,\gr}} {\gr < \gr_0}$.
\par\noindent
We have $2$ cases:
\par\noindent
$\gr_0$ is limit: Set
\begin{align*}
& \forall \gr<\gr_0 \ \gas^{\gx_0,\gr_0} > \gas^{\gx_0, \gr},
\\
& u^{\gx_0,\gr_0}= \bigunion_{\gr < \gr_0} u^{\gx_0,\gr} \union 
                \set{\ordered{\gas^{\gx_0,\gr_0},t}}
                \text{ where }\gk^0(t)=\gt_\gx.
\end{align*}
We set $T^{\gx_0,\gr_0}_0$, $T^{\gx_0,\gr_0}_1$ to anything we like. 
We do  not use them later.
\par\noindent
$\gr_0 = \gr+1$:
Let $\ordered{\gns} = s^{\gx_0,\gr}(2)$. set
\begin{align*}
& u'' = u''_1 \append u''_0 = (u^{\gx_0,\gr}_0)_{\ordered{s^{\gx_0,\gr}}},
\\
&T''_0 =   \gp^{-1}_{\gas^\gx,\gas^0} (T^0_{\ordered{\gns}}),
\\
&T''_1 = \gp^{-1}_{\mc(u''_1),\gns}       (T^0(\gns)).
\end{align*}
If there are
\begin{align*}
& q'_1 \leq^* u''_1 \union \set{T''_1},
\\
& q'_0 \leq^* u''_0 \union \set{T''_0},
\end{align*}
such that
\begin{align*}
q'_1 \append q'_0 \in D
\end{align*}
then set
\begin{align*}
\gas^{{\gx_0},\gr_0} &= \mc (q'_0),
\\
u^{{\gx_0},\gr_0}_0 &= u^{\gx_0,\gr} \union (q'_0 \setminus 
       ( u''_0 \union
                \set{T^{q'_0}} )),
\\
T^{{\gx_0},\gr_0}_0 &= T^{q'_0},
\\
u^{{\gx_0},\gr_0}_1 &= q'_1 \setminus 
        ( u''       _1 \union
                \set{T^{q'_1}} ),
\\
T^{{\gx_0},\gr_0}_1 &= T^{q'_1},
\end{align*}
otherwise set
\begin{align*}
\gas^{{\gx_0},\gr_0} &= \gas^{\gx_0, \gr},
\\
u^{{\gx_0},\gr_0}_0 &= u^{\gx_0,\gr},
\\
T^{{\gx_0},\gr_0}_0 &= T''_0,
\\
u^{{\gx_0},\gr_0}_1 &= \emptyset,
\\
T^{{\gx_0},\gr_0}_1 &= T''_1.
\end{align*}

When the induction on $\gr$ terminates we have
$\ordof{\gas^{\gx_0,\gr}, u^{\gx_0,\gr}_0, T^{\gx_0,\gr}_0,
                u^{\gx_0,\gr}_1, T^{\gx_0,\gr}_1}
                {\gr < \gt_{\gx_0}}$.
We continue with the induction on $\gx$. We set
\begin{align*}
& \forall \gr < \gt_{\gx_0}\ \gas^{\gx_0} >\gas^{\gx_0,\gr},
\\
& u^{\gx_0}_0 = \bigunion_{\gr<\gt_{\gx_0}} u^{\gx_0,\gr}_0 \union
                        \set{\ordered{\gas^{\gx_0},t}}
                        \text{ where }\max\gk^0(t)=\gt_{\gx}.
\end{align*}
When the induction on $\gx$ terminates we have
$\ordof{\gas^{\gx}, u^{\gx}_0}
                {\gx < \gk}$.
Let
\begin{align*}
& \forall \gx < \gk\ \gas^{*\prime} >\gas^{\gx},
\\
& p^{*\prime}_0 = \bigunion_{\gx<\gk} u^{\gx}_0 \union
                        \set{\ordered{\gas^{*\prime},t}}
                        \text{ where }\max\gk^0(t)=\max p_0^0.
\end{align*}
We set
\begin{align*}
& \Lev_0(T^{p^{*\prime}_0}) = \gp^{-1}_{\gas^{*\prime},\gas^{0}} \Lev_0(T^0).
\end{align*}
Let us consider $\ordered{\gns}\in \Lev_0(T^{p^{*\prime}_0})$.
There is $\gx$ such that $\gk^0(\gns) = \gt_\gx$. We set
\begin{align*}
& s(1) = \setof{\ordered{\gas,\gp_{\gas^{*\prime},\gas}(\gns)}}
		{\gas \in \supp p^*_0},
\\
& s(2)=	\set{\ordered{\gp_{\gas^{*\prime},\gas^0}(\gns)}}.
\end{align*}
Let $\gx_0=\gx+1$.
By our construction there is $\gr$ such that
\begin{align*}
& (u^{\gx_0,\gr_0}_0)_{\ordered{s}} = 
        (u^{\gx_0,\gr_0}_0)_{\ordered{s^{\gx_0,\gr}}}
\end{align*}
where $\gr_0 = \gr+1$.
We set
\begin{align*}
& T^{p^{*\prime}_0}_{\ordered{\gns}} = \gp^{-1}_{\gas^{*\prime},\gas^{\gx_0,\gr_0}}(T^{\gx_0,\gr_0})
			 \intersect
			\gp^{-1}_{\gas^{*\prime},\gas^{0}}
				(T^0_{\ordered{\gp_{\gas^*,\gas^0}(\gns)}}),
\\
& T^{p^{* \prime}_0}(\gns) =   
        \pi^{-1}_{\gns, \mc(u_1^{\gx_0,\gr_0})} T^{\gx_0,\gr_0}_1,
\\
& p'_1(\gns) =  u_1^{\gx_0,\gr_0}.
\end{align*}

Let us show that $p^{*\prime}_0$ approximates the $p^*$ 
we look for. So let $\ordered{\gns}\in T^{p^{*\prime}_0}$ and assume
\begin{subequations}\label{approx:in}
\begin{align}
& q'_1 \leq^* p'_1(\gns) \union
              ((p^{*\prime}_0)_{\ordered{\gns}})_1,
\\
& q'_0 \leq^* ((p^{*\prime}_0)_{\ordered{\gns}})_0,
\\
&q'_1 \append q'_0 \in D.
\end{align}
\end{subequations}
Let $\gx$ be such that $\gk^0(\gns) = \gt_\gx$. Set $s$ as
\begin{align*}
& s(1) = \setof{\ordered{\gas,\gp_{\gas^{*\prime},\gas}(\gns)}}
		{\gas \in \supp p^*_0},
\\
& s(2)= \set{\ordered{\gp_{\gas^{*\prime},\gas^0}(\gns)}},
\end{align*}
where $\gx_0=\gx+1$, $\gr_0=\gr+1$.
By our construction there is $\gr$ such that
\begin{align*}
(u^{\gx_0,\gr_0})_{\ordered{s^{\gx_0,\gr}}} = (u^{\gx_0,\gr_0})_{\ordered{s}}.
\end{align*}
Let us set
\begin{align*}
& r =  \big( (u^{\gx_0,\gr_0}_{\ordered{s^{\gx_0,\gr}}})_1
                        \union \set{T^{\gx_0,\gr_0}_1}
      \big)
        \append
      \big( (u^{\gx_0,\gr_0}_{\ordered{s^{\gx_0,\gr}}})_0
                \union 
                                \set{T^{\gx_0,\gr_0}_0}
        \big).
\end{align*}
By construction we have
\begin{align*}
\big( p'_1(\gns) \union ((p^{*\prime}_0)_{\ordered{\gns}})_1 \big) \append
((p^{*\prime}_0)_{\ordered{\gns}})_0
                 \leq^* r.
\end{align*}
So what we have is
\begin{align*}
D \ni q'_1 \append q'_0 \leq^* r.
\end{align*}
This is a positive answer to the question in the induction, hence
\begin{align*}
r \in D,
\end{align*}
which gives us, by openness of $D$, that
\begin{align} \label{approx:out}
\big( p'_1(\gns) \union ((p^{*\prime}_0)_{\ordered{\gns}})_1 \big) \append
((p^{*\prime}_0)_{\ordered{\gns}})_0 \in D.
\end{align}
Having proved this approximation property of $p^{*\prime}_0$,
let us consider the set
\begin{align*}
B = \setof {\ordered{\gns} \in T^{p^{*\prime}_0}}
             {\exists q \leq^* 
\big( p'_1(\gns) \union ((p^{*\prime}_0)_{\ordered{\gns}})_1 \big) \append
((p^{*\prime}_0)_{\ordered{\gns}})_0\ 
                q \in D}.
\end{align*}
Let $\gas^{*\prime}=\mc(p_0^{*\prime})$. There are 2 cases here:
\begin{enumerate}
\item $\exists \gz < \len(\Es)\ B \in E_{\gas^{*\prime}}(\gz)$.
\par\noindent
Let us set
\begin{align*}
& p^\gz_1 = j(p'_1)(\Es_{\gas^{*\prime}}\restricted \gz),
\\
& \gbs^{\gz\prime} = \mc(p^\gz_1),
\\
& A^\gz = \setof  {\ordered{\gns}}
        {\big((p^\gz_1)_{\ordered{\gns}}\big)_1 = 
                p'_1(\gp_{\gbs^{\gz\prime},\gas^{*\prime}}(\gns))}.
\end{align*}
Clearly
\begin{align*}
& A^\gz \in E_{\gbs^{\gz\prime}}(\gz).
\end{align*}
Let $\gbs^\gz$ be the larger of $\gbs^{\gz\prime}$, $\gas^{*\prime}$.
Set
\begin{align*}
& T^\gz = (\gp_{\gbs^\gz,\gas^{*\prime}})^{-1} (T^{p^{*\prime}_0}),
\\
& p^\gz = p^\gz_1 \union p^{*\prime}_0 \union \set{T^\gz}.
\end{align*}
The nice property of $p^\gz$ is that when 
$\ordered{\gns} \in T^\gz\restricted 
                        (\gp_{\gbs^\gz,\gbs^{\gz\prime}})^{-1}(A^\gz)$
we get that 
\begin{align*}
(p^\gz)_{\ordered{\gns}} \leq^* 
        \big( p'_1(\gp_{\gbs^\gz,\gbs^{\prime}}(\gns)) \union
                     ((p^{*\prime}_0)_{\ordered{\gp_{\gbs^\gz,\gas^{*\prime}}(\gns)}})_1
        \big)
\append
                     ((p^{*\prime}_0)_{\ordered{\gp_{\gbs^\gz,\gas^{*\prime}}(\gns)}})_0.
\end{align*}
We set
\begin{align*}
& p^* = p^\gz,
\\
& A = \Lev_0(T^\gz) \intersect (\gp_{\gbs^\gz,\gbs^{\gz\prime}})^{-1} A^\gz,
\end{align*}
and show that the claim is satisfied: Assume that
\begin{align*}
& \ordered{\gns} \in A,
\\
& q'_1 \leq^* ((p^*)_{\ordered{\gns}})_1,
\\
& q'_0 \leq^* ((p^*)_{\ordered{\gns}})_0,
\\
& q'_1 \append q'_0 \in D.
\end{align*}
Note that
\begin{align*}
& ((p^*)_{\ordered{\gns}})_1 \leq^* p'_1(\gp_{\gbs^\gz,\gas^{*\prime}})
         \union ((p^{*\prime}_0)_{\ordered{\gp_{\gbs^\gz,\gas^{*\prime}}(\gns)}})_1,
\\
& ((p^*)_{\ordered{\gns}})_0 \leq^*
         ((p^{*\prime}_0)_{\ordered{\gp_{\gbs^\gz,\gas^{*\prime}}(\gns)}})_0.
\end{align*}
Hence, we know that
\begin{align*}
& q'_1 \leq^* p'_1(\gp_{\gbs^\gz,\gas^{*\prime}}(\gns))
                       \union ((p^{*\prime}_0)_{\ordered{\gp_{\gbs^\gz,\gas^{*\prime}}(\gns)}})_1,
\\
& q'_0 \leq^* 
         ((p^{*\prime}_0)_{\ordered{\gp_{\gbs^\gz,\gas^{*\prime}}(\gns)}})_0,
\\
& q'_1 \append q'_0 \in D.
\end{align*}
This is the assumption \eqref{approx:in}. So from \eqref{approx:out} we know 
that
\begin{align*}
      \big( p'_1(\gp_{\gbs^\gz,\gas^{*\prime}}(\gns)) \union
      ((p^{*\prime}_0)_{\ordered{\gp_{\gbs^\gz,\gas^{*\prime}}(\gns)}})_1
      \big)
 \append
      ((p^{*\prime}_0)_{\ordered{\gp_{\gbs^\gz,\gas^{*\prime}}(\gns)}})_0
                \in D
\end{align*}
hence by openness of $D$
\begin{align*}
(p^*)_{\ordered{\gns}}
                \in D,
\end{align*}
\item $\forall \gz < \len(\Es)\ B \notin E_{\gas^{*\prime}}(\gz)$ which is the same
as saying that
\begin{align*}
\setof {\ordered{\gns} \in T^{p^{*\prime}_0}}
             {\forall q \leq^* \big(p'_1(\gns) \union 
                ((p^{*\prime}_0)_{\ordered{\gns}})_1 \big)
                \append
                ((p^{*\prime}_0)_{\ordered{\gns}})_0\ 
                q \notin D} \in \Es_{\gas^{*\prime}}.
\end{align*}
\par\noindent
In fact from the construction we can see that
\begin{align*}
\setof {\ordered{\gns} \in T^{p^{*\prime}_0}}
             {p'_1(\gns) = \emptyset}
                 \in \Es_{\gas^{*\prime}}.
\end{align*}
So we really have
\begin{align*}
A = \setof {\ordered{\gns} \in T^{p^{*\prime}_0}}
             {\forall q \leq^*
                (p^{*\prime}_0)_{\ordered{\gns}}\ 
                q \notin D} \in \Es_{\gas^{*\prime}}
\end{align*}
and the completion is quite easy now, we set
\begin{align*}
& T^{p^*} = T^{p^{*\prime}} \restricted A,
\\
& p^* = p^{*\prime}_0 \union \set{T^{p^*}}.
\end{align*}
\end{enumerate}
\end{proof}
\begin{claim} \label{canon-dense}
Let $D$ be dense open in $\PE/\Pe$, $p = p_0 \in \PE/\Pe$.
Then there is $p^* \leq^* p$
such that one and only one of the following is true:
\begin{enumerate}
\item There are $n<\gw$, $S \subseteq T^{p^*}\restricted [V_\gk]^{n}$ 
such that 
\begin{enumerate}
\item $\forall k<n\,\exists \gx<\len(\Es)\ 
                \Suc_S(\ordered{\gns_1,\dotsc,\gns_k})
                         \in E_{\mc(p^*)}(\gx)$,
\item $\forall \ordered{\gns_1,\dotsc,\gns_{n}} \in S$
                $(p^*)_{\ordered{\gns_1,\dotsc,\gns_{n}}} \in D$,
\end{enumerate}
\item $\forall n<\gw\ \forall \ordered{\gns_1,\dotsc,\gns_{n}} \in 
                        T^{p^*}
        \forall q \leq^* (p^*)_{\ordered{\gns_1,\dotsc,\gns_{n}}}\ 
               q \notin D$.
\end{enumerate}
\end{claim}
\begin{proof}
Let $p^0=p$.

Generate $p^{n+1}\leq^* p^n$ by invoking \ref{canon-dense:n} for
$n+1$ levels.

Take $\forall n<\gw\ p^* \leq p^n$.
\end{proof}
\begin{claim} \label{DenseHomogen-0}
Let $D$ be dense open in $\PE/\Pe$, $p = p_0 \in \PE/\Pe$.
Then there are $n<\gw$, $p^* \leq^* p$,
        $S \subseteq T^{p^*}\restricted [V_\gk]^{n}$ 
such that 
\begin{enumerate}
\item $\forall k<n\,\exists \gx<\len(\Es)\ 
                \Suc_S(\ordered{\gns_1,\dotsc,\gns_k})
                         \in E_{\mc(p^*)}(\gx)$,
\item $\forall \ordered{\gns_1,\dotsc,\gns_{n}} \in S$
                $(p^*)_{\ordered{\gns_1,\dotsc,\gns_{n}}} \in D$.
\end{enumerate}
\end{claim}
\begin{proof}
Towards a contradiction, let us assume that the conclusion is false.
That means that
for all $p^* \leq^* p$, for all $n < \gw$, for all
 $S \subseteq T^{p^*} \restricted [V_\gk]^n$ such that
\begin{align*}
\forall k<n \,\forall \ordered{\gns_1,\dotsc,\gns_k} \in S \,
        \exists \gx < \len(\Es)\ 
        \Suc_S(\ordered{\gns_1,\dotsc,\gns_k}) \in E_{\mc(p^*)}(\gx)
\end{align*}
we have
\begin{align*}
\exists \ordered{\gns_1,\dotsc,\gns_n} \in S \ 
                (p^*)_{\ordered{\gns_1,\dotsc,\gns_n}} \notin D.
\end{align*}
We construct a $\leq^*$-decreasing sequence as follows:
We set $p^0 = p$. We construct $p^{n+1}$ from $p^n$ using \ref{canon-dense:n}
for $n+1$ levels. Due to our assumption we get
\begin{align*}
\forall \ordered{\gns_1,\dotsc,\gns_n} \in T^{p^n} \,
                \forall q\leq^* (p^n)_{\ordered{\gns_1,\dotsc,\gns_n}} \ 
                q \notin D.
\end{align*}
We choose $p^{*\prime}$ such that $\forall n<\gw\ p^{*\prime}\leq^* p^n$ we
get
\begin{align*}
\forall n<\gw \, \forall \ordered{\gns_1,\dotsc,\gns_n} \in T^{p^{*\prime}_0} \,
                \forall q\leq^* (p^{*\prime})_{\ordered{\gns_1,\dotsc,\gns_n}} \ 
                q \notin D.
\end{align*}
Construct tree $T$ from $T^{p^{*\prime}}$ 
using \ref{skeleton}. Let us call $p^{*}$ the condition $p^{*\prime}$ with $T$ 
substituted for $T^{p^{*\prime}}$. Now if we have
\begin{align*}
q \leq p^*
\end{align*}
then there is $\ordered{\gns_1,\dotsc,\gns_n} \in T^{p^{*\prime}}$ such that
\begin{align*}
q \leq^* (p^{*\prime})_{\ordered{\gns_1,\dotsc,\gns_n}}.
\end{align*}
Hence
\begin{align*}
q \notin D.
\end{align*}
However, $D$ is dense. Contradiction.
\end{proof}
\begin{claim} \label{DenseHomogen-1}
Let $D$ be dense open in $\PE$, $p = p_1 \append p_0 \in \PE$.
Then there is $p^* \leq^* p$
such that 
\begin{multline*}
\exists S^1\ \exists n_1 \forall \ordered{\gns_{1,1},\dotsc,\gns_{1,{n_1}}} 
                \in S^1 \, \dotsc
\exists S^0\ \exists n_0 \forall \ordered{\gns_{0,1},\dotsc,\gns_{0,{n_0}}}
                \in S^0 \,
\\
         (p^*_1)_{\ordered{\gns_{1,1},\dotsc,\gns_{1,n_1}}} \append
          \dotsb \append
         (p^*_0)_{\ordered{\gns_{0,1},\dotsc,\gns_{0,n_0}}}
                         \in D
\end{multline*}
where
\begin{enumerate}
\item $S^i \subseteq T^{p^*_i} \restricted [V_\gk]^{n_i}$,
\item $\forall l<n_i \,\forall \ordered{\gns_1,\dotsc,\gns_{l}} \in S^i \,
        \exists \gx \ 
        \Suc_{S^i}(\ordered{\gns_1,\dotsc,\gns_{l}}) 
                                                \in E_{\mc(p^*_i)}(\gx)$.
\end{enumerate}
\end{claim}
\begin{proof}
Let $\ges$ be such that $p_1 \in \Pe$. 
We prove that
there are $n<\gw$, $p^*_0 \leq^* p_0$, $q_1 \leq p_1$,
        $S \subseteq T^{p^*}\restricted [V_\gk]^{n}$ 
such that 
\begin{enumerate}
\item $\forall k<n\,\exists \gx<\len(\Es)\ 
                \Suc_S(\ordered{\gns_1,\dotsc,\gns_k})
                         \in E_{\mc(p^*)}(\gx)$,
\item $\forall \ordered{\gns_1,\dotsc,\gns_{n}} \in S$
                $q_1\append (p^*)_{\ordered{\gns_1,\dotsc,\gns_{n}}} \in D$.
\end{enumerate}
Set
\begin{align*}
E = \setof {r\in \PE}{\exists q_1\ q_1 \append r \in D, q_1 \leq p_1}.
\end{align*}
This $E$ is dense open in $\PE/\Pe$: Let $r \in \PE/\Pe$. Then
$p_1 \append r \in \PE$. By density of $D$, there is $q_1 \append s \in D$
such that $q_1 \leq p_1$, $s \leq r$. By the definition of $E$, $s \in E$.
Hence $E$ is dense. Openness of $E$ is immediate from openness of $D$.

By \ref{DenseHomogen-0} there are $p^*_0 \leq^* p_0$, $S'$, $n<\gw$ such that
\begin{enumerate}
\item $\forall k<n\,\exists \gx<\len(\Es)\ 
                \Suc_{S'}(\ordered{\gns_1,\dotsc,\gns_k})
                         \in E_{\mc(p^*)}(\gx)$,
\item $\forall \ordered{\gns_1,\dotsc,\gns_{n}} \in S'$
                $(p^*_0)_{\ordered{\gns_1,\dotsc,\gns_{n}}} \in E$.
\end{enumerate}
This means that $\forall \ordered{\gns_1,\dotsc,\gns_{n}} \in S'$ there
is $q_1(\gns_1,\dotsc,\gns_n) \leq p_1$ such that
\begin{align*}
\forall \ordered{\gns_1,\dotsc,\gns_n} \in S'\ 
                q_1(\gns_1,\dotsc,\gns_n) \append
                        (p^*_0)_{\ordered{\gns_1,\dotsc,\gns_{n}}} \in D.
\end{align*}
As $\power{\Pe} < \gk$, $q_1(\gns_1,\dotsc,\gns_n)$ is in fact
almost always constant. Hence, by shrinking $S'$ to $S$ and letting $q_1$ be this
constant value, we get
\begin{align*}
\forall \ordered{\gns_1,\dotsc,\gns_n} \in S\ 
                q_1 \append
                        (p^*_0)_{\ordered{\gns_1,\dotsc,\gns_{n}}} \in D.
\end{align*}
With this, we finished the first part of the proof. We use this claim
for all conditions in $\Pe$.

Let $\Pe = \setof{p_1^\gz}{\gz<\gl}$ where $\gl < \gk$.

We construct by induction a $\leq^*$-decreasing sequence
$\ordof{p_0^\gz}{\gz < \gl}$.
Set
\begin{align*}
p^0_0 = p_0.
\end{align*}
Assume we have constructed $\ordof{p^\gz_0}{\gz<\gz_0}$.
\par\noindent
$\gz_0$ is limit: Choose $p^{\gz_0}_0 \leq^* p^\gz_0$ for all $\gz<\gz_0$.
\par\noindent
$\gz_0=\gz+1$: Use the first part of the proof on $p^\gz_1 \append p^\gz_0$ 
to construct $p^{\gz_0}_0$.

When the induction terminates we have $\ordof{p_0^\gz}{\gz < \gl}$.
Choose
\begin{align*}
\forall \gz<\gl\ p^*_0 \leq^* p^\gz_0.
\end{align*}
Let
\begin{align*}
D_{\ges} = \setof{q_1 \in \Pe}{\exists n\ \exists S\ 
                \forall \ordered{\gns_1,\dotsc,\gns_n} \in S\ 
                q_1 \append (p^*_0)_{\ordered{\gns_1,\dotsc,\gns_n}} \in D}.
\end{align*}
$D_{\ges}$ is dense open: Let $q_1 \in \Pe$. Then there is $\gz$ such that
$q_1=p_1^\gz$. By the induction we have that there are $n$, $S$, $r_1 \leq q_1$ 
such that
\begin{align*}
\forall \ordered{\gns_1,\dotsc,\gns_n} \in S\ 
          r_1 \append (p^{\gz_0}_0)_{\ordered{\gns_1,\dotsc,\gns_n}} \in D.
\end{align*}
By openness of $D$ we get
\begin{align*}
\forall \ordered{\gns_1,\dotsc,\gns_n} \in S\ 
           r_1 \append (p^*_0)_{\ordered{\gns_1,\dotsc,\gns_n}} \in D.
\end{align*}
Hence
\begin{align*}
           r_1 \in D_{\ges}.
\end{align*}
As $D_{\ges}$ is dense open we can use \ref{DenseHomogen-0}. Hence
there are $p^*_1 \leq p_1$, $S^1$, $n_1$ such that
\begin{align*}
\forall \ordered{\gns_1,\dotsc,\gns_{n_1}} \in S^1\ 
           (p^*_1)_{\ordered{\gns_1,\dotsc,\gns_{n_1}}} \in D_{\ges}.
\end{align*}
This means that
\begin{multline*}
\forall \ordered{\gns_{1,1},\dotsc,\gns_{1,n_1}} \in S^1\ 
\exists S\ \exists n\ 
\forall \ordered{\gns_{0,1},\dotsc,\gns_{0,n}} \in S\ 
\\
(p^*_1)_{\ordered{\gns_{1,1},\dotsc,\gns_{1,n_1}}} \append
(p^*_0)_{\ordered{\gns_{0,1},\dotsc,\gns_{0,n}}} \in D,
\end{multline*}
which is what we need to prove.
\end{proof}
Finally, we add the last touch.
\begin{proof}[Proof of \ref{DenseHomogen}]
The proof is done by induction on $k$. The case $k=1$ is \ref{DenseHomogen-1}.
We assume, then, that the theorem is proved for $k$ and prove it for $k+1$.

Then $p = p_{k+1}\append p_k \append \dotsb p_0$.
Let $\ges$ be such that $p_{k+1} \in \Pe$. 
We just repeat the proof of \ref{DenseHomogen-1} with $\Pe$ and use
the induction hypothese to conclude the proof.
\end{proof}
%
% Prikry Condtion Section
%
\section{Prikry's condition} \label{Prikry'sCondition}
\begin{theorem} \label{PrikryCondition}
Let $p \in \PE$ and $\gs$ a formula in the forcing language. Then there
is $p^* \leq p$ such that $p^* \decides \gs$.
\end{theorem}
\begin{proof}
The set $\setof{q \in \PE}{q \decides \gs}$ is dense open. 
Assuming $p = p_k \append \dotsb \append p_0$ and using \ref{DenseHomogen}
we get that there is $q \leq^* p$ such that
\begin{multline*}
\exists S^k\ \exists n_k \forall \ordered{\gns_{k,1},\dotsc,\gns_{k,{n_k}}} 
                \in S^k \, \dotsc
\exists S^0\ \exists n_0 \forall \ordered{\gns_{0,1},\dotsc,\gns_{0,{n_0}}}
                \in S^0 \,
\\
         (p^*_k)_{\ordered{\gns_{k,1},\dotsc,\gns_{k,n_k}}} \append
          \dotsb \append
         (p^*_0)_{\ordered{\gns_{0,1},\dotsc,\gns_{0,n_0}}}
                         \decides \gs.
\end{multline*}
Recall that we really should write 
\begin{align*}
& S^{k-1}(\gns_{k,1},\dotsc,\gns_{k,n_k}),
\\
& S^{k-2}(\gns_{k,1},\dotsc,\gns_{k,n_k},\gns_{k-1,1},\dotsc,\gns_{k-1,n_k}),
\\
& \vdots
\end{align*}
In order to avoid (too much) clutter, we use the following convention
in the proof. When we write 
\begin{align*}
        \gnv \in \prod_{1\leq l\leq k} S^l
\end{align*}
we mean that
\begin{align*}
&\ordered{\gns_{k,1},\dotsc,\gns_{k,n_k}} \in S^k,
\\
&\vdots
\\
&\ordered{\gns_{1,1},\dotsc,\gns_{1,n_1}} \in S^1,
\end{align*}
and $r(\gnv)$ is
\begin{align*}
(q_k)_{\ordered{\gns_{k,1},\dotsc,\gns_{k,n_k}}}\append\dotsb\append
    (q_1)_{\ordered{\gns_{1,1},\dotsc,\gns_{1,n_1}}}.
\end{align*}
We start by naming $q_0$ as $q^{n_0}_0$ and $T^{q_0}$ as $T^{0,n_0}$.
For $\ordered{\gns_1,\dotsc,\gns_{n_0-1}} \in S^0$ set
\begin{align*}
& A^0=
\setof{\ordered{\gns_{n_0}} \in 
                \Suc_{T^{0,n_0}}(\ordered{\gns_{1},\dotsc,\gns_{n-1}})}
      {r(\gnv)\append (q^{n_0}_0)_{\ordered{\gns_{1},\dotsc,\gns_{n_0}}} \forces \gs},
\\
& A^1=
\setof{\ordered{\gns_{n_0}} \in 
        \Suc_{T^{0,n_0}}(\ordered{\gns_{1},\dotsc,\gns_{n_0-1}})}
      {r(\gnv)\append (q^{n_0}_0)_{\ordered{\gns_{1},\dotsc,\gns_{n_0}}} \forces \lnot\gs}.
\end{align*}
Note that
\begin{align*}
& \Suc_{S^0}(\ordered{\gns_{1},\dotsc,\gns_{n_0-1}})
        \subseteq
        A^0 \union
        A^1,
\\
& A^0\intersect A^1 = \emptyset
\end{align*}
Hence, there is $\gx < \len(\Es)$ such that one and only one of the following
is true:
\begin{enumerate}
\item $A^0
                        \in E_{\mc(q_0^{n_0})}(\gx)$,
\item $A^1
                        \in E_{\mc(q_0^{n_0})}(\gx)$.
\end{enumerate}
In either case, using \ref{fill-missing}, we can shrink
$T^{0,n_0}_{\ordered{\gns_1,\dotsc,\gns_{n_0-1}}}$
and get a condition $q'_0$ such that
$r(\gnv)\append q'_{\ordered{\gns_1,\dotsc,\gns_{n_0-1}}} \decides \gs$.

So we shrink now $T^{0,n_0}_{\ordered{\gns_1,\dotsc,\gns_{n_0-1}}}$ for all
$\ordered{\gns_1,\dotsc,\gns_{n_0-1}} \in S^0$ and we call this 
tree $T^{0,n_0-1}$.
The name of the condition $q^{n_0}_0$ with $T^{0,n_0-1}$ substituted for 
$T^{0,n_0}$ is $q^{n_0-1}_0$.

$q^{n_0-1}_0$ satisfies
\begin{align*}
\forall \ordered{\gns_{1},\dotsc,\gns_{n_0-1}} \in S^0 \ 
     r(\gnv)\append (q^{n_0-1}_0)_{\ordered{\gns_1,\dotsc,\gns_{n_0-1}}} 
                \decides \gs.
\end{align*}
We are now in the same position as we 
were when setting $q^{n_0}_0$. So by
repeating the above arguments we get
\begin{align*}
p_0 \ge^* q_0=q^{n_0}_0 \ge^* q^{n_0-1}_0 \ge^* \dotsb \ge^* q^1_0 \ge^* q^0_0
\end{align*}
such that for each $l = n_0,n_0-1,\dotsc,1,0$
\begin{align*}
\forall \ordered{\gns_1,\dotsc,\gns_l} \in S^0 \ 
           r(\gnv)\append (q^l_0)_{\ordered{\gns_1,\dotsc,\gns_l}} \decides \gs.
\end{align*}
Specifically we get
\begin{align*}
           r(\gnv)\append (q^0_0) \decides \gs.
\end{align*}
Of course $q^0_0$ depends on $\gnv$.
Note that we got from $q^{n_0}_0$ to $q^0_0$ only by shrinking the trees.
So, we repeat this process for all $\gnv$
calling the resulting condition $q^0_0(\gnv)$.
So we have
\begin{align*}
\forall \gnv\in\prod_{1\leq l\leq k}S^l\ 
        r(\gnv)\append q^0_0(\gnv) \decides \gs.
\end{align*}
By setting
\begin{align*}
T^{p^*_0}=\bigintersect_{\gnv} T^{q^0_0(\gnv)}
\end{align*}
and letting $p^*_0$ be $q_0$ with $T^{p^*_0}$ substituted for $T^{q_0}$ we get
\begin{align*}
\forall \gnv\in\prod_{1\leq l\leq k}S^l\ 
        r(\gnv) \append p^*_0 \decides \gs.
\end{align*}

We are in the same position as in the beginning of the proof. So we can
generate in the same way $p^*_1$ from $p_1$ and so on until we have
\begin{align*}
        p^*_k\append\dotsb\append p^*_0 \decides \gs.
\end{align*}
\end{proof}
%
% W e a k    p r o p e r n e s s
%
%
%
%
%
%
\section{Properness} \label{Properness}
The notions $\ordered{N,P}$-generic and properness are due
to Shelah \cite{Shelah}.
\begin{definition}
Let $N \subelem H_\gc$ such that
\begin{enumerate}
\item $\power{N}=\gk$,
\item $N \supseteq V_\gk$,
\item $N \supseteq N^{\upto\gk}$,
\item $P \in N$.
\end{enumerate}
Then $p \in P$ is called $\ordered{N,P}$-generic if
\begin{align*}
   p \forces \formula{\forall D\in \VN{N}\ D \text{ is dense open in }\VN{P}\ 
        \implies D \intersect \GN{G} \intersect 
                \VN{N} \not= \emptyset}.
\end{align*}
\end{definition}
\begin{definition}
A forcing notion $P$ is called proper if for all
$N \subelem H_\gc$ such that
\begin{enumerate}
\item $\power{N}=\gk$,
\item $N \supseteq V_\gk$,
\item $N \supseteq N^{\upto\gk}$,
\item $P \in N$,
\end{enumerate}
and for all $q \in P \intersect N$ there is $p \leq q$ which
is $\ordered{N, P}$-generic.
\end{definition}
\begin{claim} \label{NPgeneric-1}
Let $p \in \PE$, $N \subelem H_\gc$ such that
\begin{enumerate}
\item $\power{N}=\gk$,
\item $N \supseteq V_\gk$,
\item $N \supseteq N^{\upto\gk}$,
\item $\PE \in N$,
\item $p \in \PE \intersect N$.
\end{enumerate}
Then there is $p^* \leq^* p$ such that $p^*$ is $\ordered{N, \PE}$-generic
\end{claim}
\begin{proof}
Let $p = p_{k(p)}\append \dotsb \append p_1 \append p_0$.

Let $\ordof{D_\gx}{\gx < \gk}$ be an enumeration of all dense open subsets
of $\PE$ which are in $N$. Note that for $\gx_0 < \gk$ we have that
$\ordof{D_\gx}{\gx < \gx_0} \in N$.

We start now an induction on $\gx$ in which we build
\begin{align*}
\ordof{\gas^\gx, u^\gx}{\gx < \gk}.
\end{align*}
The construction is done ensuring that
$\ordof{\gas^\gx, u^\gx}{\gx < \gx_0} \in N$ for all $\gx_0 < \gk$.
We start by setting
\begin{align*}
u^0 &= p_0 \setminus \set{T^{p_0}},
\\
\gas^0 &= \mc(p_0),
\\
T^0 &= T^{p_0} \restricted \pi^{-1}_{\gas^0,0}\setof{\gns}
                                     {\gk^0(\gns) \text{ is inaccessible}},
\end{align*}
and taking an increasing enumeration in $N$
\begin{align*}
\setof{\gk^0(\gns)}{\ordered{\gns} \in T^0}
                = \ordof{\gt_\gx}{\gx < \gk}.
\end{align*}
Assume then that we have
\begin{align*}
\ordof{\gas^\gx, u^\gx}{\gx < \gx_0}.
\end{align*}
The constructions splits now according to wether $\gx_0$ is limit
or successor. In both cases the work is done inside $N$.
\par\noindent
$\gx_0$ is limit: Choose $\gas^{\gx_0} > \gas^\gx$ for all $\gx < \gx_0$
and set
\begin{align*}
u^{\gx_0} &= \bigunion_{\gx < \gx_0} u^\gx \union
        \set{\ordered{\gas^{\gx_0}, t}}
                \text{ where } \gk^0(t) = \gt_{\gx_0}.
\end{align*}
\par\noindent
$\gx_0 = \gx + 1:$
For each $\gns_1,\dotsc,\gns_n$ such that
$\gk^0(\gns_1)<\dotsb<\gk^0(\gns_n)=\gt_\gx$ we set
\begin{align*}
 S(\gns_1,\dotsc,\gns_n) = &
        \big( \prod_{\gas \in \supp u^\gx} 
                \setof {\gms_1}{\gk^0(\gms_1)=\gk^0(\gns_1)} \big)
                                     \times
\\
        & \big( \prod_{\gas \in \supp u^\gx} 
                \setof {\gms_2}{\gk^0(\gms_2)=\gk^0(\gns_2)} \big)
                                        \times 
\\
        & \vdots
\\
        & \big( \prod_{\gas \in \supp u^\gx} 
                \setof {\gms_n}{\gk^0(\gms_n)=\gk^0(\gns_n)} \big)
                                        \times 
\\
        & \set {\ordered{\gns_1,\dotsc,\gns_n}}.
\end{align*}
Let
\begin{align*}
S = \bigunion_{\gk^0(\gns_1)<\dotsb<\gk^0(\gns_n)=\gt_\gx} 
        S(\gns_1,\dotsc,\gns_n) 
\end{align*}
and set enumeration of $S$
\begin{align*}
S=\ordof{s^{\gx_0,\gr}} {\gr < \gt_{\gx_0}}.
\end{align*}
We do induction on $\gr$ which  builds
\begin{align*}
\ordof{\gas^{\gx_0,\gr}, u_0^{\gx_0,\gr},T_0^{\gx_0,\gr}}
      {\gr < \gt_{\gx_0}},
\end{align*}
from which we build $\ordered{\gas^{\gx_0}, u^{\gx_0}}$.
Set
\begin{align*}
& \gas^{\gx_0,0}=\gas^\gx,
\\
& u^{\gx_0,0}_0 = u^\gx_0.
\end{align*}
Assume we have constructed
$\ordof{\gas^{\gx_0,\gr}, u_0^{\gx_0,\gr},T_0^{\gx_0,\gr}} {\gr < \gr_0}$.
\par\noindent
$\gr_0$ is limit: Set
\begin{align*}
& \forall \gr<\gr_0 \ \gas^{\gx_0,\gr_0} > \gas^{\gx_0, \gr},
\\
& u^{\gx_0,\gr_0}= \bigunion_{\gr < \gr_0} u^{\gx_0,\gr} \union 
                \set{\ordered{\gas^{\gx_0,\gr_0},t}}
                \text{ where }\gk^0(t)=\gt_{\gx_0}.
\end{align*}
We set $T^{\gx_0,\gr_0}$ to anything we like as we do not use it later.
\par\noindent
$\gr_0 = \gr+1$:
Let $s^{\gx_0,\gr}(n(s)+1)=\ordered{\gns_1,\dotsc,\gns_n}$.
\begin{align*}
& u'' = (u^{\gx_0,\gr}_0)_{\ordered{s^{\gx_0,\gr}}},
\\
&T''_0 =   \gp^{-1}_{\gas^\gx,\gas^0} (T^0_{\ordered{\gns_1,\dotsc,\gns_n}}),
\\
&T''_1 = \gp^{-1}_{\mc(u''_1),\gns_n}
        (T^0_{\ordered{\gns_1,\dotsc,\gns_{n-1}}}(\gns_n)),
\\
& \vdots
\\
&T''_{n-1} = \gp^{-1}_{\mc(u''_{n-1}),\gns_2} (T^0_{\ordered{\gns_1}}(\gns_2)),
\\
&T''_n = \gp^{-1}_{\mc(u''_n),\gns_1} (T^0(\gns_1)).
\end{align*}
Take enumeration
\begin{multline*}
        \setof{D_\gs}{\gs< \gt_{\gx}} \times
\\
 \setof{q}{q \leq p_{k(p)}\append\dotsb\append p_1\append
                u''_{n(s)}\union \set{T''_{n}}\append\dotsb
                \append u''_1 \union \set{T''_1}
           } =
\\
        \ordof{ \ordered{ E^{\gx_0,\gr_0,\gz}, q^{\gx_0,\gr_0,\gz} } }
                { \gz<\gt_{\gx_0} }.
\end{multline*}
\par\noindent
We start induction on $\gz$. Set
\begin{align*}
& \gas^{\gx_0,\gr_0,0} = \gas^{\gx_0,\gr},
\\
& u^{\gx_0,\gr_0,0}_0 = u^{\gx_0,\gr}.
\end{align*}
Assume we have constructed
$\ordof{\gas^{\gx_0,\gr_0,\gz},u^{\gx_0,\gr_0,\gz}_0,T^{\gx_0,\gr_0,\gz}_0}
                {\gz<\gz_0}$.
\par\noindent
$\gz_0$ is limit:
\begin{align*}
& \forall \gz<\gz_0\ \gas^{\gx_0,\gr_0,\gz_0} > \gas^{\gx_0,\gr_0,\gz} ,
\\
& u^{\gx_0,\gr_0,\gz_0}_0 = \bigunion_{\gz< \gz_0} u^{\gx_0,\gr_0,\gz}_0 \union
        \set{\ordered{\gas^{\gx_0,\gr_0,\gz_0},t}}\text{ where }
                        \gk^0(t)=\gt_{\gx_0}.
\end{align*}
We set $T^{\gx_0,\gr_0,\gz}$ to whatever we want as no use of it is made
later.
\par\noindent
$\gz_0=\gz+1$:
We set
\begin{align*}
& u'' = (u_0^{\gx_0,\gr_0,\gz})_{s^{\gx_0,\gr}},
\\
& T''_0 = \gp^{-1}_{\gas^{\gx_0,\gr_0,\gz},\gas^0} T^0_{\ordered{\gns_1,
				\dotsc,\gns_n}}.
\end{align*}
If there is
\begin{align*}
& u'_0 \leq^* u''_0 \union \set{T''_0}
\end{align*}
such that
\begin{align*}
q^{\gx_0,\gr_0,\gz} \union u'_0 \in E^{\gx_0,\gr_0,\gz}
\end{align*}
then set
\begin{align*}
\gas^{\gx_0,\gr_0,\gz_0} &= \mc (u'_0),
\\
u^{\gx_0,\gr_0,\gz_0}_0 &= u^{\gx_0,\gr_0,\gz}_0 \union (u'_0 \setminus 
       ( u''_0 \union
                \set{T^{u'_0}} )),
\\
T^{\gx_0,\gr_0,\gz_0} &= T^{u'_0},
\end{align*}
otherwise set
\begin{align*}
\gas^{\gx_0,\gr_0,\gz_0} &= \gas^{\gx_0,\gr_0,\gz},
\\
u^{\gx_0,\gr_0,\gz_0}_0 &= u^{\gx_0,\gr_0,\gz}_0,
\\
T^{\gx_0,\gr_0,\gz_0}_0 &= T''_0.
\end{align*}
When the induction on $\gz$ terminates we have
$\ordof{\gas^{\gx_0,\gr_0,\gz}, u^{\gx_0,\gr_0,\gz}_0, T^{\gx_0,\gr_0,\gz}}
                {\gz < \gt_{\gx_0}}$

We continue with the induction on $\gr$. We set
\begin{align*}
& \forall \gz < \gt_{\gx_0}\ \gas^{\gx_0,\gr_0} >\gas^{\gx_0,\gr_0,\gz} ,
\\
& u^{\gx_0,\gr_0}_0 = \bigunion_{\gz<\gt_{\gx_0}} u^{\gx_0,\gr_0,\gz}_0 \union
                        \set{\ordered{\gas^{\gx_0,\gr_0},t}}
                        \text{ where }\gk^0(t)=\gt_{\gx_0}.
\end{align*}
When the induction on $\gr$ terminates we have
$\ordof{\gas^{\gx_0,\gr}, u^{\gx_0,\gr}_0, T^{\gx_0,\gr}}
                {\gr < \gt_{\gx_0}}$.
We continue with the induction on $\gx$. We set
\begin{align*}
& \forall \gr < \gt_{\gx_0}\ \gas^{\gx_0} >\gas^{\gx_0,\gr} ,
\\
& u^{\gx_0}_0 = \bigunion_{\gr<\gt_{\gx_0}} u^{\gx_0,\gr}_0 \union
                        \set{\ordered{\gas^{\gx_0},t}}
                        \text{ where }\gk^0(t)=\gt_{\gx_0}.
\end{align*}
When the induction on $\gx$ terminates we have
$\ordof{\gas^{\gx}, u^{\gx}_0}
                {\gx < \gk}$.
We note that this sequence is not in $N$. 
Let
\begin{align*}
& \forall \gx < \gk\ \gas^* >\gas^{\gx} ,
\\
& p^*_0 = \bigunion_{\gx<\gk} u^{\gx}_0 \union
                        \set{\ordered{\gas^{*},t}}
                        \text{ where } \gk^0(t)=\max p_0^0.
\end{align*}
We construct a series of trees, $R^n$, and
$T^{p^*_0}$ is $\bigintersect_{n<\gw} R^n$.
\begin{align*}
& \Lev_0(R^0) = \gp^{-1}_{\gas^{*},\gas^{0}} \Lev_0(T^0).
\end{align*}
Let us consider $\ordered{\gns_1}\in \Lev_0(R^0)$.
There is $\gx$ such that $\gk^0(\gns_1) = \gt_\gx$. We set
\begin{align*}
& s(0) = \setof{\ordered{\gas,\gp_{\gas^*,\gas}(\gns_1)}}
		{\gas \in \supp p^*_0},
\\
& s(1)=	\set{\ordered{\gp_{\gas^*,\gas^0}(\gns_1)}}.
\end{align*}
Let $\gx_0=\gx+1$.
By our construction there is $\gr$ such that
\begin{align*}
& (u^{\gx_0,\gr_0}_0)_{\ordered{s}} = 
        (u^{\gx_0,\gr_0}_0)_{\ordered{s^{\gx_0,\gr}}},
\end{align*}
where $\gr_0 = \gr+1$.
We set
\begin{align*}
& R^1_{\ordered{\gns_1}} = \gp^{-1}_{\gas^*,\gas^{\gx_0,\gr_0}}(T^{\gx_0,\gr_0})
			 \intersect
			\gp^{-1}_{\gas^*,\gas^{0}}
				(T^0_{\ordered{\gp_{\gas^*,\gas^0}(\gns_1)}}).
\end{align*}
Assume that we have constructed $R^n$.
We set the first $n$ levels of $R^{n+1}$ to be the same as the first
$n$ levels of $R^n$ and we complete the tree as follows.
Let us consider $\ordered{\gns_1,\dotsc,\gns_n}\in R^n$.
There is $\gx$ such that $\gk^0(\gns_n) = \gt_\gx$. We set $s$ as folows
\begin{align*}
& \forall 1\leq k\leq n\ 
	s(k) = \setof{\ordered{\gas,\gp_{\gas^*,\gas}(\gns_k)}}
		{\gas \in \supp p^*_0},
\\
&  s(n+1)= \set{\ordered{\gp_{\gas^*,\gas^0}(\gns_1),\dotsc,
			\gp_{\gas^*,\gas^0}(\gns_n)}}.
\end{align*}
Let $\gx_0=\gx+1$.
By our construction there is $\gr$ such that
\begin{align*}
& (u^{\gx_0,\gr_0}_0)_{\ordered{s}} = (u^{\gx_0,\gr_0}_0)_
			{\ordered{s^{\gx_0,\gr}}},
\end{align*}
where $\gr_0 = \gr+1$.
We set
\begin{align*}
& R^{n+1}_{\ordered{\gns_1,\dotsc,\gns_n}} = 
		\gp^{-1}_{\gas^*,\gas^{\gx_0,\gr_0}}(T^{\gx_0,\gr_0})
			 \intersect
			\gp^{-1}_{\gas^*,\gas^{0}}
		(T^0_{\ordered{\gp_{\gas^*,\gas^0}(\gns_1),\dotsc,
			\gp_{\gas^*,\gas^0}(\gns_n)}}).
\end{align*}
After $\gw$ stages we set
\begin{align*}
& T^{p^*_0} = \bigintersect_{n<\gw} R^n.
\end{align*}
We finish the construction by setting
\begin{align*}
& p^* = p_{k(p)} \append \dotsb \append p_1 \append p^*_0.
\end{align*}
We show that $p^*$ is as required.

Let $G$ be $\PE$-generic such that $p^* \in \PE$. Let $D\in N$ be
dense open in $\PE$. We want to show that $D\intersect G\intersect N
\not= \emptyset$.

Choose $q \append r_0 \in D \intersect G$ such that
\begin{align*}
& r_0 \leq^* p''_0,
\\
& q \leq p_{k(p)} \append \dotsb \append p_1 \append
		p''_{n} \append \dotsb \append p''_1,
\end{align*}
where
\begin{align*}
& \ordered{\gns_1,\dotsc,\gns_n} \in \dom T^{p^*_0},
\\
& \gk^0(\gns_n) = \gt_\gx,
\\
& D \in \ordof{D_\gz}{\gz<\gt_\gx},
\\
& p'' = (p^*_0)_{\ordered{\gns_1,\dotsc,\gns_n}}.
\\
\end{align*}
We set $s$ to be
\begin{align*}
& \forall 1\leq k\leq n\ 
	s(k) = \setof{\ordered{\gas,\gp_{\gas^*,\gas}(\gns_k)}}
		{\gas \in \supp p^*_0},
\\
&  s(n+1)= \set{\ordered{\gp_{\gas^*,\gas^0}(\gns_1),\dotsc,
			\gp_{\gas^*,\gas^0}(\gns_n)}}.
\end{align*}
We get that
\begin{align*}
        (p^*_0)_{\ordered{\gns_1,\dotsc,\gns_n}} = (p^*_0)_{\ordered{s}}.
\end{align*}
We let $\gx_0 = \gx+1$.
Recall the enumeration of $S$ in the construction. There is $\gr$
such that
\begin{align*}
        (p^*_0\setminus \set{\ordered{\mc(p^*_0),(p^*_0)^{\mc}}})_{\ordered{s}}
	\union \set{\ordered{\mc(p^*_0),(p^*_0)^{\mc}}} = 
		(p^*_0)_{\ordered{s^{\gx_0,\gr}}}.
\end{align*}
We let $\gr_0 = \gr +1$.
Considering the construction of $T^{p^*_0}$, we see that
\begin{align*}
& T^{p^*_0}_{\ordered{\gns_1,\dotsc,\gns_n}}\leq T^{\gx_0,\gr_0}
\end{align*}
hence 
\begin{align*}
 (p^*_0)_{\ordered{s}} \leq^* (u^{\gx_0,\gr_0}_0)_{\ordered{s^{\gx_0,\gr_0}}}.
\end{align*}							    
We note that
\begin{align*}
\forall 1\leq k\leq n\ p''_{k} = ((p^*_0)_{\ordered{s}})_k.
\end{align*}							    
Recalling that $q$ was chosen so that
\begin{align*}
q \leq p_{k(p)} \append \dotsb \append p_1 \append 
		p''_{n} \append \dotsb \append p''_1
\end{align*}							    
we conclude that there is $\gz$ such that $q = q^{\gx_0,\gr_0,\gz}$ 
and $D=E^{\gx_0,\gr_0,\gz}$.
That is
\begin{align*}
E^{\gx_0,\gr_0,\gz} \ni q^{\gx_0,\gr_0,\gz} \append r_0 \leq^* 
		q^{\gx_0,\gr_0,\gz} \append ((u^{\gx_0,\gr_0,\gz})_
			{\ordered{s^{\gx_0,\gr}}})_0.
\end{align*}							    
Note that this is an answer to the question we asked in the construction.
Hence, due to elementarity of $N$, there was such a condition in $N$.
Hence
\begin{align*}
q^{\gx_0,\gr_0,\gz} \append 
	((u^{\gx_0,\gr_0,\gz})_{\ordered{s^{\gx_0,\gr}}})_0 \in
	D \intersect N.
\end{align*}							    
The last point to note is that
\begin{align*}
q^{\gx_0,\gr_0,\gz} \append 
	((u^{\gx_0,\gr_0,\gz})_{\ordered{s^{\gx_0,\gr}}})_0 \geq^*.
q^{\gx_0,\gr_0,\gz} \append p''_0 \geq^* q\append r_0 \in G
\end{align*}							    
Hence
\begin{align*}
q^{\gx_0,\gr_0,\gz} \append 
	((u^{\gx_0,\gr_0,\gz})_{\ordered{s^{\gx_0,\gr}}})_0 \in G.
\end{align*}							    
\end{proof}
\begin{corollary} \label{properness}
$\PE$ is proper.
\end{corollary}
\section{Cardinals in $V^{\PE}$} \label{Cardinals}
\begin{lemma} \label{nocollapse-special}
$\gk^+$ remains a cardinal in $V^{\PE}$.
\end{lemma}
\begin{proof}
The proof really has no connection to the specific structure of $P_\Es$. It is
an exercise in properness.

Let
\begin{align*}
p \forces \formula{\GN{f} \func \VN{\gk} \to \VN{\gk^+}}.
\end{align*}
Choose $\gc$ large enough so that $H_\gc$ contains everything we are
interested in. Take $N \subelem H_\gc$ such that
\begin{enumerate}
\item $p, \PE, \GN{f} \in N$,
\item $\power{N}= \gk$,
\item $N \supseteq V_\gk$,
\item $N \supseteq N^{\upto\gk}$.
\end{enumerate}
By \ref{properness} there is $q \leq p$ which is $\ordered{N,\PE}$-generic.
Let us set
\begin{align*}
\gl = N \intersect \gk^+,
\end{align*}
where $\gl$ is an ordinal $< \gk^+$.

Let $G$ be $\PE$-generic with $q\in G$.
The $\ordered{N,\PE}$-genericity ensures us that for all $\gx < \gk$ 
        $\GN{f}(\gx)^{N[G]} \in N$ and
        $\GN{f}(\gx)^{N[G]}= \GN{f}(\gx)^{H_\gc[G]}$.
Hence $\ran \GN{f}^{V[G]} \subseteq \gl$. That is
\begin{align*}
q \forces \formula{\GN{f} \text{ is bounded in } \gk^+}.
\end{align*}
\end{proof}
\begin{lemma} \label{nocollapse-above}
No cardinals $> \gk$ are collapsed by $\PE$.
\end{lemma}
\begin{proof}
$\gk^+$ is not collapsed by \ref{nocollapse-special}.
No cardinals $\geq \gk^{++}$ are collapsed as 
$\PE$ satisfies $\gk^{++}$-c.c.
\end{proof}
\begin{lemma}
\label{small-forcing}
Let $\gx < \gk$ and $\gz$ the ordinal such that
$\gk^0(\Es_G(\gz)) \leq \gx < \linebreak[4] \gk^0(\Es_G(\gz + 1))$. Then
$\Pset(\gx) \intersect V[G] = \Pset(\gx) \intersect V[G\restricted \gz]$.
\end{lemma}
\begin{proof}
Take $p = p_n \append \dotsb \append p_{k+1} \append p_k \append \dotsb \append
p_0 \in G$ such that $\Es(p_{k+1}) = \Es_G(\gz)$,
$\Es(p_{k}) = \Es_G(\gz + 1)$.
We know that $V[G] = V[G/p]$.
So we work in $P_{\Es}/p$. Set $p^l = p_n \append \dotsb \append p_{k+1}$,
$p^h = \set{\ordered{\Es_G(\gz), \emptyset}} \append p_k \append \dotsb
        \append p_0$. Then $P_{\Es}/p = P_{\Es} \restricted p^l \times
        P_{\Es} \restricted p^h$.
Note that $\ordered{P_{\Es}/p^h, \leq^*}$ is $\gk^0(\Es_G(\gz+1))$-closed.
In particular it is $\gx^+$-closed.

Let $A \in V[G]$, $A \subseteq \gx$.
Choose $\GN{A}$, a canonical $P_{\Es}/p$-name for $A$.
Let $q \in P_{\Es}/p^h$. By induction we construct 
$\ordof{q_\gt}{\gt < \gx}$ satisfying
\begin{enumerate}
\item $\gt_0 < \gt_1 \implies q_{\gt_1} \leq^* q_{\gt_0}$,
\item $q_\gt \decides \formula{ \VN{\gt} \in \GN{A}}$.
\end{enumerate}
Choose $q_\gx \leq^* q_\gt$ for all $\gt < \gx$.

By density argument we can construct, $\GN{B}$, a $P\restricted p^l$-name
such that $A = \GN{A}[G] = \GN{B}[G\restricted p^l]$.
\end{proof}
\begin{corollary}
No carindals $\leq \gk$ are collapsed by $\PE$.
\end{corollary}
\begin{proof}
Let $G \subseteq \PE$ be generic. Assume $\gl < \gk$ is a collapsed cardinal.
Let $\gm = \power{\gl}^{V[G]}$. We have $\gm < \gl$, and there is
$A \in \Pset(\gm)^{V[G]}$ which codifies the order type $\gl$.
Let $\gz$ be the
unique ordinal such that $\gk^0(\Es_G(\gz)) \leq \gm < \gk^0(\Es_G(\gz+1))$.
By \ref{small-forcing} $A \in V[G\restricted \gz]$. Hence $\gl$ is collapsed
already in $V[G\restricted \gz]$. However, by \ref{nocollapse-above},
$P_{\Es_G(\gz)}$ collapses no cardinals above $\gk^0(\Es_G(\gz))$.
Contradiction.

So, no cardinal $<\gk$ is collapsed. As $\gk$ is a limit of cardinals
which are not collapsed, it is not collapsed.
\end{proof}
We have just shown
\begin{theorem}
No cardinals are collapsed in $V^{\PE}$.
\end{theorem}
\section{Properties of $\gk$ in $V^{P_{\Es}}$} \label{PropertiesK}
\begin{theorem}
If $\len(\Es) = \gk^+$ then $V^{P_{\Es}} \satisfies \formula{ \gk 
\text{ is regular} }$.
\end{theorem}
\begin{proof}
Let $\gl < \gk$, $\GN{f}$ be such that
\begin{align*}
\forces_{P_\Es} \formula{ \GN{f}\func \VN{\gl} \to \VN{\gk}}.
\end{align*}
Let
\begin{align*}
D_0 = \setof{p} {\exists i \ p\forces \formula{ \GN{f}(0)=\VN{i} }}.
\end{align*}
As $D_0$ is a dense open set we can invoke \ref{DenseHomogen} to get
$p^{\prime 0}$, $n_0$, $S^0 \subset T^{p^{\prime 0}}$ such that
\begin{align*}
&\forall k<n_0 \,\forall \ordered{\gn_1,\dotsc,\gn_k} \in S^{\prime 0} \,
        \exists \gx < \len(\Es)\ 
        \Suc_S(\ordered{\gn_1,\dotsc,\gn_k}) \in E_{\mc(p^{\prime 0})}(\gx),
\\
&\forall \ordered{\gn_1,\dotsc,\gn_{n_0}} \in S^{\prime 0} \ 
                (p^{\prime 0})_{\ordered{\gn_1,\dotsc,\gn_{n_0}}} \in D_0.
\end{align*}
Let us set
\begin{align*}
A'_0 = \setof { (p^{\prime 0})_{\ordered{\gn_1,\dotsc,\gn_{n_0}}} }
             { \ordered{\gn_1,\dotsc,\gn_{n_0}} \in S^{\prime 0}}.
\end{align*}
$A_0$ is an anti-chain. By shrinking $T^{p^{\prime 0}}$ as was done in
the proof of \ref{PrikryCondition} we can make $A_0$ into a 
maximal anti-chain below $p^{\prime 0}$.
As $\gl < \gk$ and $\ordered{P^*_{\Es},\leq^*}$ is $\gk$-closed
we can construct a $\leq^*$-decreasing sequence
\begin{align*}
p^{\prime 0} \geq^* p^{\prime 1} \geq^* \dotsb \geq^* 
                p^{\prime \gt} \geq^* \dotsb \qquad \gt < \gl.
\end{align*}
and $n_\gt$, $S^{\prime \gt} \subseteq T^{p^{\prime \gt}}$ such that
\begin{align*}
&\forall k<n_\gt \,\forall \ordered{\gn_1,\dotsc,\gn_k} \in S^{\prime \gt} \,
        \exists \gx < \len(\Es)\ 
        \Suc_{S^{\prime \gt}}(\ordered{\gn_1,\dotsc,\gn_k}) \in E_{\mc(p^{\prime \gt})}(\gx),
\\
&\forall \ordered{\gn_1,\dotsc,\gn_{n_\gt}} \in S^{\prime \gt} \ 
  \exists i \ 
    (p^{\prime \gt})_{\ordered{\gn_1,\dotsc,\gn_{n_\gt}}} \forces
        \formula { \GN{f}(\VN{\gt})=\VN{i} },
\end{align*}
and 
\begin{align*}
A'_\gt = \setof { (p^{\prime \gt})_{\ordered{\gn_1,\dotsc,\gn_{n_\gt}}} }
             { \ordered{\gn_1,\dotsc,\gn_{n_\gt}} \in S^{\prime \gt}}
\end{align*}
is a maximal anti-chain below $p^{\prime \gt}$.

Let $p' \leq^* p^{\prime \gt}$ for all $\gt < \gl$.
We set 
$S^\gt = \gp^{-1}_{\mc(p'),\mc(p^{\prime\gt)}}(S^{\prime\gt})$ and 
$p^\gt$ to be $p'$ with
$\gp^{-1}_{\mc(p'),\mc(p^{\prime\gt)}}(T^{p^{\prime\gt}})$ substituted for 
$T^{p'}$ and maybe shrunken a bit so that
\begin{align*}
A_\gt = \setof { (p^{\gt})_{\ordered{\gn_1,\dotsc,\gn_{n_\gt}}} }
             { \ordered{\gn_1,\dotsc,\gn_{n_\gt}} \in S^{\gt}}
\end{align*}
is a maximal anti-chain below $p^\gt$.

Let $p \leq^* p^\gt$ for all $\gt < \gl$ and let $\GN{g}$ be the following
$\PE$-name
\begin{align*}
\GN{g} = \bigunion_{\gt<\gl}
            \setof { \ordered{\VN{\ordered{\gt, i}}, 
                                (p^\gt)_{\ordered{\gn_1,\dotsc,\gn_{n_\gt}}}} }
                { A_\gt \ni (p^\gt)_{\ordered{\gn_1,\dotsc,\gn_{n_\gt}}}\forces 
                        \formula{\GN{f}(\gt)= \VN{i}}\ 
                    }.
\end{align*}
Then
\begin{align*}
p\forces \formula { \GN{f}=\GN{g} }.
\end{align*}
Let $P^*$ be the following forcing notion:
\begin{align*}
P^* = \setof {q \leq_{R} p} {q \in \PE}.
\end{align*}
By \ref{SubForcing} $\ordered{P^*, \leq_R}$ is sub-forcing of 
$\ordered{\PE/p, \leq}$.
Hence if $G$ is $\PE$-generic then
$G^* = G \intersect P^*$ is $P^*$-generic. $\GN{g}$ is in fact a
$P^*$-name and as can be seen from its' definition $\GN{g}[G] =
\GN{g}[G^*] \in V[G^*]$. So in order to complete the proof it is
enough to show that $\forces_{P^*} \formula{\GN{g} \text{ is bounded}}$.

By \ref{GenericForRadin} there is $r \in R_{\mc(p)}$
such that $P^* \simeq R_{\mc(p)}/r$. Now we use the following
fact about Radin forcing: When the measure sequence is of length
$\gk^+$, $\gk$ is regular in the generic extension. Necessarily,
$\forces_{P^*} \formula{
\GN{g} \text{ is bounded}}$.
\end{proof}
\begin{definition}
We say that $\gt < \len(\Es)$ is a repeat point of $\Es$ if
$P_\Es=P_{\Es\restricted\gt}$.
\end{definition}
Note that if $\gt$ is a repeat point then $P_{\Es\restricted\gt}\in M$.
\begin{theorem}
If $\Es$ has a repeat point, $j''\setof{\GN{A}_\gx}{\gx<\gl} \in M$
where $\linebreak[0]\setof{\GN{A}_\gx}{\gx<\gl}$ is an enumeration of all canonical
$\PE$-names of subsets of $\gk$, then 
$V^{\PE} \satisfies \linebreak[0] 
        \ulcorner \gk \text{ is}\linebreak[0] \text{ measurable} \urcorner$.
\end{theorem}
\begin{proof}
We use the usual method under these circumstances. Let
$\gt$ be a repeat point of $\Es$ and
$G$ be $\PE$-generic over $V$. For the duration of this proof
let us define: 
\begin{itemize}
\item If $p=p_0 \in P^*_{\Es}$ then
        \begin{align*}
        p\restricted\gt = \setof{\ordered{\Es_\ga\restricted\gt,p^{\Es_\ga}}}
                                {\Es_\ga \in \supp p} \union \set{T^p}.
        \end{align*}
\item
        If $p=p_n\append\dotsb\append p_1\append p_0$ then
        \begin{align*}
        p\restricted\gt = p_n\append\dotsb\append p_1\append p_0\restricted\gt.
        \end{align*}
\end{itemize}
Let us set
\begin{align*}
G\restricted\gt = \setof{p\restricted\gt}{p \in G}.
\end{align*}
We note that
\begin{enumerate}
\item $\gl < j(\gk)$,
\item $M\satisfies\formula {j(P^*_{\Es})/\PE\text{ is }j(k)-closed}$,
\item $G\restricted\gt$ is $P_{\Es\restricted\gt}$-generic over $M$.
\end{enumerate}
So, we can construct a $\leq^*$-decreasing sequence in $M[G\restricted\gt]$,
\ordof{p^\gx}{\gx<\gl}, such that
\begin{align*}
p^{\gx+1} \decides \formula{\gk \in j(\GN{A}_\gx)}.
\end{align*}
We define an ultrafilter, $U$, in $V$ by
\begin{align*}
\GN{A}_\gx[G] \in U \iff p_{\gx+1}\forces \formula {\gk \in j(\GN{A}_\gx)}.
\end{align*}
\end{proof}
The assumptions we used in this theorem are very strong. We believe that
a repeat point is enough in order to get measurability.
\ifnum\resurrect = 1
Let $\gt$ be a repeat point of $\Es$. Hence $\PE = P_{\Es \restricted \gt}$.
For the duration of this proof let us use the following definition:
For $p \in P^*_{\Es}$
\begin{align*}
p\restricted \gt = \setof {\ordered{\Es_\ga\restricted \gt, p^{\Es_\ga}}} 
                          {\Es_\ga \in \supp p}
\end{align*}
Note that $P_{\Es\restricted\gt} \in M$ and if $p \in \PE$ then
$p\restricted \gt \union \set{T^p} \in P_{\Es\restricted \gt}$.

Let $G$ be $\PE$-generic over $V$. By setting
\begin{align*}
G\restricted \gt = \setof {p\restricted \gt \union \set{T^p}} 
                          {p \in G}
\end{align*}
we get that $G\restricted \gt$ is $P_{\Es\restricted\gt}$-generic over $M$.

Let $p\in P^*_{\Es}$. Then $j(p) \in j(P^*_{\Es})$. Due to the definition 
of extender
sequences we get that ${\mc(p)\restricted\gt} \in j(T^p)$. Hence
$j(p)_{\ordered{{\mc(p)\restricted\gt}}} \in j(\PE)$.
Specifically we note that if we write
\begin{align*}
j(p)_{\ordered{{\mc(p)\restricted\gt}}} = p_1 \union p_0
\end{align*}
then $p_1 = p\restricted \gt \union \set {j(T^p)(\mc(p)\restricted\gt)}$.
Moreover, ${j(T^p)(\mc(p)\restricted\gt)}$ is an $\mc(p)$-tree, as
$\gt$ is a repeat point.

In what follows we need to handle differently the part of a condition
which is moved by $j$ and the parts which are not moved. So, when
writing $q\append p$ we mean that $q\append p = p_n\append\dotsb\append
p_1\append p_0$ where $p = p_0$ and 
$q = p_n\append\dotsb\append p_1$.
We define $U$ in $V[G]$ as follows:
\begin{align*}
\GN{A}[G] \in U \iff \exists q\append p\in G \ & \exists S
\\
              &  q\append j(p)_{\ordered{\mc(p)\restricted\gt}}
                        \setminus \set {j(T^p)} \union \set{S}
                 \forces
                        \formula{\VN{\gk} \in \GN{A}}.
\end{align*}
We show that $U$ is a measure.
\begin{enumerate}
\item The definition makes sense: Let $q^1\append p^1,q^2\append p^2 \in G$,
        $S^1, S^2$ trees. Let $q\append p \in G$ such that
        $q\append p \leq q^1\append p^1,q^2\append p^2$. Then
\begin{align*}
q\append j(p)_{\ordered{\mc(p)\restricted\gt}}
                    \setminus \set {j(T^{p})} \union& \set{S^1 \intersect S^2}
        \leq
\\
&q^1\append j(p^1)_{\ordered{\mc(p^1)\restricted\gt}}
                        \setminus \set {j(T^{p^1})} \union \set{S^1},
\\
&q^2\append j(p^2)_{\ordered{\mc(p^2)\restricted\gt}}
                        \setminus \set {j(T^{p^2})} \union \set{S^2}
\end{align*}
\item Completeness: Let $\GN{A}^{\gx}[G] \in U$ for all $\gx<\gk$.
By our definition there are $q^\gx\append p^\gx \in G$, $S^\gx$ such that
\begin{align*}
q^\gx\append j(p^\gx)_{\ordered{\mc(p^\gx)\restricted\gt}}
                        \setminus \set {j(T^{p^\gx})} \union \set{S^\gx}
        \forces \formula{\VN{\gk} \in j(\GN{A}^\gx)}
\end{align*}
The set
\begin{align*}
D = \setof {q\append p \in \PE} {\forall \gx<\gk\  
                                q\append p \le q^\gx\append p^\gx \text{ or }
              \exists \gx<\gk\  q\append p \incompatible q^\gx \append p^\gx}
\end{align*}
is dense, hence there is $q\append p \in G$ such that $\forall \gx<\gk\ 
q\append p \leq q^\gx\append p^\gx$. Let
\begin{align*}
S = \bigintersect_{\gx < \gk} S^\gx
\end{align*}
As $M \supseteq M^\gk$ we have $S\in M$. We got that for all $\gx<\gk$
\begin{align*}
q\append j(p)_{\ordered{\mc(p)\restricted\gt}}
                        \setminus \set {j(T^{p})} \union& \set{S}
\leq
& q^\gx\append j(p^\gx)_{\ordered{\mc(p^\gx)\restricted\gt}}
                        \setminus \set {j(T^{p^\gx})} \union \set{S^\gx}
\end{align*}
Hence
\begin{align*}
q\append j(p)_{\ordered{\mc(p)\restricted\gt}}
                        \setminus \set {j(T^{p})} \union& \set{S}
 \forces\formula{\forall \gx<\gk\ \VN{\gk}\in j(\GN{A}^\gx)}
\end{align*}
That is
\begin{align*}
\dintersect_{\gx<\gk} \GN{A}^\gx[G] \in U
\end{align*}
\end{enumerate}
\fi
\section{What have we proved?} \label{TheTheorem}
We can sum everything as follows:

We can control independently two properties of $\gk$ in a generic extension.
The first is the size of $2^\gk$ which is controlled by $\power{\Es}$.
The second is how `big' we want $\gk$ to be which is controlled by
$\len(\Es)$.
%
%
%
%
%
% Generic by Iteration
%
\section{Generic by Iteration} \label{ByIteration}
Recall that if $R$ is Radin forcing generated from $j\func V \to M$ then
there is $\gt$ and $G \in V$ such that $G$ is $j_{0,\gt}(R)$-generic
over $M_\gt$.

Our original aim was to find some form of this claim for our forcing.
We have a partial result in this direction. Namely, when $\len(\Es)=1$
we have a generic filter in $V$ over an elementary submodel in $M_\gw$.

In this section we assume that $\len(\Es)=1$.

Let us take an iteration of $j=j_{0,1}$
\begin{equation*}
\ordered{
        \ordof{M_n}{n < \gw},
        \ordof{j_{n,m}}{n \leq m < \gw}
        }.
\end{equation*}
Choose $\gc$ large enough so that everything interesting
is in $H_\gc$ (i.e. $P_{\Es} \in H_\gc$)
and set
\begin{align*}
& \gk_n = j_{0,n}(\gk),
\\
& \Es^n = j_{0,n}(\Es),
\\
& P^n = j_{0,n}(P_\Es)(=P_{\Es^n}),
\\
& \gc_n = j_{0,n}(\gc).
\end{align*}
\begin{definition}[when $\len(\Es)=1$]
We call \ordered{N, p} a $k$-pair if
\begin{enumerate}
\item $M_k \satisfies \formula{ N \subelem H_{\gc_k}^{M_k}}$,
\item $M_k \satisfies \formula{ \power{N} = \gk_k }$,
\item $M_k \satisfies \formula{ N \supseteq V_{\gk_k}^{M_k} }$,
\item $M_k \satisfies \formula{ N \supseteq N^{\upto \gk_k} }$,
\item $p \in P^{k+1} \intersect j_{k,k+1}(N) $,
\item If $D \in N$, $D$ is dense open in $P^k$
        then there are 
        $n$, $S\leq T^p$ such that
	\begin{align*}
          \forall \ordered{\gn_1,\dotsc,\gn_n}\in S \ 
                (p)_{\ordered{\gn_1,\dotsc,\gn_n}} \in j_{k,k+1}(D).
	\end{align*}
\end{enumerate}
\end{definition}
\begin{claim} \label{find-pair-1}
Let $N \in M_k$, $k<\gw$ and $q \in j_{k,k+1}(N) \intersect P^{k+1}$ 
such that $M_k$ satisifes:
\begin{enumerate}
\item $N \subelem H_{\gc_k}^{M_k}$,
\item $\power{N} = \gk_k$,
\item $N \supseteq V^{M_k}_{\gk_k}$,
\item $N \supseteq N^{\upto \gk_k}$,
\item $P^k \in N$.
\end{enumerate}
Then there is
        $p \in j_{k,k+1}(N) \intersect P^{k+1}$ such that
\begin{enumerate}
\item $p \leq^* q$,
\item $\ordered{N,p}$ is a $k$-pair.
\end{enumerate}
\end{claim}
\begin{proof}
Let \ordof{D^k_{\gx}}{\gx < \gk_k} be an enumeration, in $M_k$, of all the
dense open subsets of $P^k$ which are in $N$. Set $N_{k+1} = j_{k,k+1}(N)$.
As we have
\begin{align*}
& \ordof{j_{k,k+1}(D^k_{\gx})}{\gx < \gk_k} \subseteq N_{k+1},
\\
& \ordof{j_{k,k+1}(D^k_{\gx})}{\gx < \gk_k} \in M_{k+1},
\\
&M_{k+1} \satisfies \formula{ N_{k+1} \supseteq N_{k+1}^{\upto \gk_{k+1}}
        },
\end{align*}
we get that 
\begin{align*}
& \ordof{j_{k,k+1}(D^k_{\gx})}{\gx < \gk_k} \in N_{k+1}.
\end{align*}
Starting with $q$ we construct in $N_{k+1}$ a $\leq^*$-decreasing sequence
$\ordof{p_\gx}{\gx<\gk_k}$ using \ref{DenseHomogen}. Note
that we have no problem at the limit stages as
$N_{k+1} \satisfies \formula{ \ordered{P^{k+1},\leq^*} \text{ is } \gk_{k+1}
\text{-closed} }$. Choose now $p \in P^{k+1} \intersect N_{k+1}$
such that $\forall \gx < \gk_{k}\ p \leq^* p_\gx$. We get that
$\ordered{N,p}$ is a $k$-pair.
\end{proof}
%
% Generic Approximation Sequence
%
\begin{definition}
We call $\Cond{N}{p}$  a $P^\gw$-generic approximation sequence if
\begin{align*}
\Cond{N}{p} = \ordof{\ordered{N_k, p^{k+1}}}{k_0 \leq k < \gw}
\end{align*}
such that for all $k_0 \leq k < \gw$
\begin{enumerate}
\item $M_{k} \satisfies \formula{ N_k \subelem H_{\gc_{k}}^{M_{k}}}$,
\item $M_{k} \satisfies \formula{ \power{N_k} = \gk_{k} }$,
\item $M_{k} \satisfies \formula{ N_k \supseteq V^{M_{k}}_{\gk_{k}} }$,
\item $M_{k} \satisfies \formula{ N_k \supseteq N^{\upto \gk_{k}} }$,
\item $P^{k} \in N_k$,
\item $\ordered{N_k,p^{k+1}}$ is a $k$-pair,
\item $j_{k_1, k_2}(N_{k_1}) = N_{k_2}$,
\item $p^{k+2} \leq^* (j_{k+1,k+2}(p^{k+1}))_{\ordered{\mc(p^{k+1})}}$.
\end{enumerate}
\end{definition}
\begin{definition}
Let $\Cond{N}{r}$ be a $P^\gw$-generic approximating sequence. Then
\begin{align*}
G(\Cond{N}{r})=\setof{p \in P^\gw}{\exists k\ j_{k,\gw}(r^k) \leq^* p}.
\end{align*}
\end{definition}
%
% Given N's we can complete to a GenSeq
% GenSeqComplete-1 exists
%
\begin{claim} \label{GenSeqComplete-1}
Let ${k_0} < \gw$, $q \in P^{k_0}\intersect N_{k_0}$
and assume for all $k_0\le k < \gw$
\begin{enumerate}
\item $P^{k} \in N_{k}$,
\item $j_{k,k+1}(N_k) = N_{k+1}$,
\item $M_{k} \satisfies \formula{ N_{k} \subelem H_{\gc_{k}}^{M_{k}}}$,
\item $M_{k} \satisfies \formula{ \power{N_{k}} = \gk_{k} }$,
\item $M_{k} \satisfies \formula{ N_{k} \supseteq V^{M_{k}}_{\gk_{k}} }$,
\item $M_{k} \satisfies \formula{ N_{k} \supseteq N_{k}^{\upto \gk_{k}} }$.
\end{enumerate}
Then there is a $P^\gw$-generic approximating sequence
\begin{align*}
\Cond{N}{p} = \ordof{\ordered{N_k, p^{k+1}}}{k_0 \leq k < \gw}
\end{align*}
such that $p^{{k_0}+1} \leq^* j_{k_0,k_0+1}(q)$.
\end{claim}
\begin{proof}
We construct the $p^{k+1}$ by induction. We set $p^{k_0} = q$.

Assume that we have constructed $\ordof{p^{k'}}{k_0 \leq k'\leq k}$.
We set $\linebreak[0] q^{k+1} = \linebreak[4] j_{k,k+1}(p^k)_{\ordered{\mc(p^k)}}$.
Invoke \ref{find-pair-1} to get $p^{k+1} \leq^* q^{k+1}$ such that
$\ordered{N_k, p^{k+1}}$ is a $k$-pair.

When the induction terminates we have
         $\ordof{\ordered{N_k, p^{k+1}}}{k_0 \leq k < \gw}$ as required.
\end{proof}
\begin{claim} \label{GenSeq-1}
Let ${k_0} < \gw$ and assume
\begin{enumerate}
\item $M_{k_0} \satisfies \formula{ N_{k_0} \subelem H_{\gc_{k_0}}^{M_{k_0}}}$,
\item $M_{k_0} \satisfies \formula{ \power{N_{k_0}} = \gk_{k_0} }$,
\item $M_{k_0} \satisfies \formula{ N_{k_0} \supseteq V^{M_{k_0}}_{\gk_{k_0}} }$,
\item $M_{k_0} \satisfies \formula{ N_{k_0} \supseteq N_{k_0}^{\upto \gk_{k_0}} }$,
\item $P^{{k_0}} \in N_{k_0}$,
\item $q \in P^{k_0}\intersect N_{k_0}$.
\end{enumerate}
Then there is a $P^\gw$-generic approximating sequence
\begin{align*}
\Cond{N}{p} = \ordof{\ordered{N_k, p^{k+1}}}{k_0 \leq k < \gw}
\end{align*}
such that $p^{{k_0}+1} \leq^* j_{k_0,k_0+1}(q)$.
\end{claim}
\begin{proof}
We set $N_k = j_{k_0, k}(N_{k_0})$ for all $k_0 < k < \gw$ and then
we invoke \ref{GenSeqComplete-1}.
\end{proof}
%
%       N-generic 1
%
\begin{theorem} \label{N-generic-1}
Assume
\begin{enumerate}
\item $M_\gw \satisfies \formula{ N \subelem H_{\gc_\gw}^{M_\gw}}$,
\item $M_\gw \satisfies \formula{ \power{N} = \gk_\gw }$,
\item $M_\gw \satisfies \formula{ N \supseteq V^{M_\gw}_{\gk_\gw} }$,
\item $M_\gw \satisfies \formula{ N \supseteq N^{\upto \gk_\gw} }$,
\item $P^{\gw} \in N$,
\item $q \in P^{\gw}\intersect N$.
\end{enumerate}
Then there is, in $V$, a filter $G \subseteq P^\gw$ 
such that 
\begin{enumerate}
\item $q \in G$,
\item $\forall D \in N$ D is dense open in $P^\gw$ $\implies$
        $G \intersect D \intersect N \not= \emptyset$.
\end{enumerate}
\end{theorem}
\begin{proof}
We find, $\Cond{N}{p}$, a $P^\gw$-generic approximating sequence.
$G(\Cond{N}{p})$ is the required filter.

Find $k_0$ and $N_{k_0}$, $q^{k_0}$ such that
\begin{align*}
& P^{k_0} \in N_{k_0},
\\
& j_{k_0,\gw}(N_{k_0}) = N,
\\
& j_{{k_0},\gw}(q^{k_0}) = q.
\end{align*}
Invoke \ref{GenSeq-1} to get from $N^{k_0}$, $q^{k_0}$ a generic
approximating sequence $\Cond{N}{p}$.

Let $D \in N$ be dense open in $P^\gw$. 

Find $k \geq k_0$ and $D^k$ such that $j_{k,\gw}(D^k)=D$.  
As usual we set
$D^{k+l} = j_{k,k+l}(D^k)$.
By construction there is $n$
such that
\begin{align*}
\forall \ordered{\gn_1,\dotsc,\gn_n}\in T^{p^{k+1}}\,
        (p^{k+1})_{\ordered{\gn_1,\dotsc,\gn_n}} \in D^{k+1} \intersect N_{k+1},
\end{align*}
which means that
\begin{align*}
j_{k+1,k+1+n}(p^{k+1})_{\ordered{\mc(p^{k+1}),\dotsc,j_{k+1,k+n}(\mc(p^{k+1}))}} 
        \in D^{k+1+n} \intersect N_{k+1+n}.
\end{align*}
Hence
\begin{align}\label{GoverN:dense}
j_{k+1,\gw}(p^{k+1})_{\ordered{\mc(p^{k+1}),\dotsc,j_{k+1,k+n}(\mc(p^{k+1}))}} 
        \in D \intersect N.
\end{align}
As $\Cond{N}{p}$ is a $P^\gw$-generic approximation sequence it satisfies
\begin{align*}
p^{k+1+n} \leq^* 
j_{k+1,k+1+n}(p^{k+1})_{\ordered{\mc(p^{k+1}),\dotsc,j_{k+1,k+1+n}(\mc(p^{k+1}))}}.
\end{align*}
Hence
\begin{multline*}
j_{k+1+n,\gw}(p^{k+1+n}) \leq^*
\\
j_{k+1,\gw}(p^{k+1})_{\ordered{\mc(p^{k+1}),\dotsc,j_{k+1,k+n}(\mc(p^{k+1}))}},
\end{multline*}
giving us that
\begin{align}\label{GoverN:generic}
j_{k+1,\gw}(p^{k+1})_{\ordered{\mc(p^{k+1}),\dotsc,j_{k+1,k+n}(\mc(p^{k+1}))}}
        \in G(\Cond{N}{p}),
\end{align}
so $G({\Cond{N}{p}}) \intersect D \intersect N \not= \emptyset$ by \ref{GoverN:dense}
and \ref{GoverN:generic}.
\end{proof}
%
%
%  Concluding Remarks
%
%
\section{Concluding Remarks} \label{ConcludingRemarks}
\begin{enumerate}
\item
We believe a repeat point is enough in order to get measurability
of $\gk$. We use a much stronger assumption in our proof.
\item
A definition of repeat point that depends only on the
extender sequence and is equivalent to the one we gave
(which mentions $\PE$) will probably be useful.
\item
It is not completly clear what $\len(\Es)$ should be in order
to make sure that $\Es$ has a repeat point.
\item
A finer analysis in the case of measurability and stronger properties
is needed. For example, extending the elementary embedding to the
generic extension, and not just constructing a normal ultrafilter.
\item
We do not know how to get a generic by iteration when
$\len(\Es)>1$.
\item
Making this forcing more `precise' by adding `gentle' collapses
so we get a prescribed behaviour on \emph{all} cardinals below
$\gk$ in the generic extension is in preparation.
\end{enumerate}
%
%

%\section{preliminaries}
%\begin{rem} Begining remark \end{rem}
%\begin{definition} 1st definition \end{definition}
%\begin{theorem} interesting theorem \end{theorem}
%\begin{lemma} interesting lemma \end{lemma}
%\begin{lemma} 2nd interesting lemma \end{lemma}
%\begin{theorem} 2nd interesting theorem \end{theorem}
%\begin{definition} 2nd definition \end{definition}
%\begin{rem} Final remark \end{rem}
\end{document}